\documentclass[12pt]{amsart}
\usepackage{amsmath,amscd,amssymb,amsfonts,latexsym}
\usepackage[all]{xy}
\usepackage{ulem}

\usepackage[latin1]{inputenc}

\usepackage[T1]{fontenc}

\usepackage{eucal} 
\usepackage{cmmib57} 
\usepackage{amsbsy}
\usepackage{amsthm,alltt}

\usepackage{graphicx,enumerate}

\usepackage{geometry}
\newcommand{\ST}{\mathcal{S}}
\newcommand{\CP}{\mathbb{CP}^{n}}
\newcommand{\CN}{\mathbb{C}^{n}}
\newcommand{\U}{\mathcal{U}}
\newcommand{\C}{\mathbb{C}}
\newcommand{\Z}{\mathbb{Z}}
\newcommand{\Q}{\mathbb{Q}}
\newcommand{\F}{\mathcal{F}}
\newcommand{\D}{\mathcal{D}}
\newcommand{\K}{\mathcal{L}}
\newcommand{\R}{\mathcal{\Re}^{\bullet}}

\newcommand{\g}{\mathbf{g}}
\newcommand{\f}{\mathbf{f}}
\newcommand{\V}{\mathcal{V}}
\newcommand{\W}{\mathcal{W}}
\newcommand{\fst}{\textbf{1}}
\newcommand{\Su}{{\rm Supp}}
\newcommand{\Hom}{{\rm Hom}}
\newcommand{\Ext}{{\rm Ext}}
\newcommand{\T}{\mathbb{T}}

\newcommand{\re}{\mathcal{R}}

\newtheorem{thm}{Theorem}[section]
\newtheorem{prop}[thm]{Proposition}
\newtheorem{lem}[thm]{Lemma}
\newtheorem{cor}[thm]{Corollary}
\theoremstyle{definition}
\newtheorem{definition}[thm]{Definition}
\newtheorem{example}[thm]{Example}
\theoremstyle{remark}
\newtheorem{remark}[thm]{Remark}
\numberwithin{equation}{section}

\def\be{\begin{equation}}
\def\ee{\end{equation}}

\def\bt{\begin{thm}}
\def\et{\end{thm}}

\def\bc{\begin{cor}}
\def\ec{\end{cor}}

\def\br{\begin{remark}}
\def\er{\end{remark}}

\def\bp{\begin{prop}}
\def\ep{\end{prop}}

\def\bl{\begin{lem}}
\def\el{\end{lem}}

\def\bex{\begin{example}}
\def\eex{\end{example}}

\def\bd{\begin{definition}}
\def\ed{\end{definition}}

\title{Characteristic Varieties of Hypersurface Complements}
\author{Yongqiang Liu}
\address{Y. Liu: School of Mathematical Sciences, University of Science and Technology of China, No.96,  JinZhai Road Baohe District, Hefei, Anhui, 230026, P.R.China}
\email {liuyq@mail.ustc.edu.cn}

\author{Lauren\c{t}iu Maxim}
\address{L. Maxim: Department of Mathematics, University of Wisconsin, 480 Lincoln Drive, Madison, WI 53706, USA}
\email {maxim@math.wisc.edu}

\date{\today}

\keywords{free abelian cover, hypersurface complement, hyperplane arrangements, non-isolated singularities,  Alexander varieties, characteristic varieties, Dwyer-Fried sets, resonance varieties.}

\subjclass[2010]{14J17, 14J70, 32S20, 32S25, 32S55}

\begin{document}

\begin{abstract}   We give divisibility results for the (global) characteristic varieties of hypersurface complements expressed in terms of the local characteristic varieties at points along one of the irreducible components of the hypersurface.  As an application, we recast old and obtain new finiteness and divisibility results for the classical (infinite cyclic) Alexander modules of complex hypersurface complements.  Moreover, for the special case of hyperplane arrangements, we translate our divisibility results for characteristic varieties in terms of the corresponding resonance varieties.
 \end{abstract}

\maketitle

\tableofcontents

\section{Introduction}
\subsection{Foreword}
The study of the topology of complex algebraic varieties is a classical subject going back to Zariski, Enriques, Deligne, Hirzebruch or Arnol'd.  A particularly fruitful approach consists of using local topological information at singular points in order to extract global topological properties of singular algebraic varieties.

A celebrated result in this circle of ideas is Libgober's divisibility theorem for the classical (infinite cyclic) Alexander invariants of complements to complex affine hypersurfaces with only isolated singularities, including at infinity (see \cite{L82, L83, L94}). More precisely, Libgober showed that the (only nontrivial) global  Alexander polynomial of such a  complement divides the product of the local Alexander polynomials associated with each singular point.  The second-named author  (\cite{Max}) extended these divisibility results to the case of hypersurfaces with arbitrary singularities and in general position at infinity. Furthermore, this global-to-local approach was used in \cite{DM} to show that such divisibility results also hold for certain multivariable Alexander-type invariants, called the {Alexander varieties} (or {support loci}) of the hypersurface complement.

The Alexander varieties of a topological space are closely related to the so-called 
characteristic varieties (or jumping loci of rank-one local systems), which are topological invariants of the space.  A fundamental result of Arapura (\cite{A}) showed that the characteristic varieties of complex hypersurface complements have a rigid structure. More precisely, the characteristic varieties of a plane curve complement are unions
of subtori of the character torus, possibly translated by unitary characters. The result is  known to be true more generally, for any smooth quasi-projective variety, by work of Budur-Wang  \cite{BW}. Furthermore, Libgober \cite{L09} and Budur-Wang \cite{BW15} extended this rigidity property also in the local context of a hypersurface singularity germ complement.  The interplay between the characteristic  and  Alexander varieties of any topological space is investigated in detail in \cite{PS10}.

\subsection{Main results}
Let $f:\C^{n+1} \to \C$ ($n\geq 2$) be a reduced homogeneous polynomial of degree $d>1$. Then $f$ defines a global Milnor fibration $f: M\to \C^{\ast}$ with total space $$M=\C^{n+1} \setminus f^{-1}(0)$$ and Milnor fiber $F=f^{-1}(1)$.  Let $V$ be the hypersurface in $\CP$ defined by $f$, and set $$M^{\ast}=\CP \setminus V.$$  Assume that $f=f_{1}\cdot \cdots \cdot f_{r}$, where $f_{i}$ are the irreducible factors of $f $ with degree $d_{i}$. 
Then $V_{i}=\lbrace f_{i}=0 \rbrace$ are the irreducible components of $V$. 

Let $\V_{i}(M)$ and $\V_{i}(M^{\ast})$ denote the $i$-th characteristic varieties of $M$ and $M^{\ast}$, respectively (see Definition \ref{d2.3}). These are sub-varieties of the $r$-dimensional torus $(\C^{\ast})^{r}$, which is the scheme of the multivariable Laurent polynomial ring $\Gamma_{r}:=\C[t_{1},t_{1}^{-1},\cdots,t_{r},t_{r}^{-1}]$.

\medskip

The aim of this paper is to establish a relationship between the characteristic varieties of $M$ and $M^\ast$, and to compute them in terms of the corresponding local invariants around singularities via divisibility-type results.

\medskip

One of our main results, Theorem \ref{t5.1}, relates the characteristic varieties of $M$ and $M^\ast$ as follows:
\bt  For any $k\geq 0$, \be 
\bigcup_{i=0}^{k} \V_{i}(M)= \bigcup_{i=0}^{k} \V_{i}(M^{\ast})\subset V(\prod_{i=1}^{r} t_{i}^{d_{i}}=1).
\ee
\et
\noindent (Here, for an ideal $I$ of a commutative ring $R$, we denote as usual by $V(I)$  the reduced subscheme of $Spec({R})$ defined by $I$.)

\medskip

We also derive general divisibility-type results expressing the {\it global} characteristic varieties $\V_{i}(M)$ and $\V_{i}(M^{\ast})$ of $M$ and $M^{\ast}$, respectively, in terms of the {\it local} topology around singularities. 
Let us describe some of these results in more detail. 
Let $x \in V$ be a point, and $B_{x}$ be a small open ball in $\CP$ centered at $x$. Set $M^{\ast}_{x}=M^{\ast}\cap B_{x}$ for the {\it local} complement. Set $I_{x}=\lbrace i \mid f_{i}(x)= 0 \rbrace$. 
In order to account for the possibility that some $f_{i}$ does not vanish at $x$ (i.e., $i\notin I_{x}$), we introduce the {\it uniform (local) Alexander varieties} $\W_i^{unif}(M^{\ast}_{x}, \nu_{x})$, see Definition \ref{d2.1} and Definition  \ref{d3.0} for details. We denote by $\W(M^{\ast}_{x}, \nu_{x})$, $\W^{unif}(M^{\ast}_{x}, \nu_{x})$  and $\V(M^{\ast}_{x})$ the union of the corresponding (uniform) Alexander varieties and characteristic varieties, respectively. 
Then one of our main divisibility results can be stated as follows (see Theorem \ref{t5.2} for a precise formulation):

\bt \label{t1.2}   With the above notations, for any $k\leq n-1$ we have that:  \be  
 \bigcup_{i=0}^{k}  \V_{i}(M)= \bigcup_{i=0}^{k}  \V_{i}(M^{\ast})  \subset  \bigcup_{S\subset V_{1}} \W^{unif}(M^{\ast}_{x_{S}}, \nu_{x_{S}}) ,
\ee
where the union is over the strata of $V_{1}$ with $\dim_{\C} S\geq n-k-1$ and $x_{S}\in S$ is arbitrary.
\et
When coupled with results of  \cite{PS10} about the relation between local Alexander and characteristic varieties (cf. Theorem 2.16), Theorem \ref{t1.2} above implies that  the characteristic varieties (resp. Alexander varieties) of $M$ and $M^{\ast}$ are (up to degree $n-1$) contained in the union of the local characteristic varieties at points along one of the irreducible components of $V$.

Our divisibility results are optimal in the sense that no additional assumptions on the hypersurace are needed.  Moreover, our results improve on  the divisibility theorems for the Alexander varieties obtained by Dimca and the second-named author in \cite{DM}, see Remarks \ref{r4.2} and \ref{r5.1}.

\medskip

In the special case of hyperplane arrangements, the above divisibility results can be sharpened as follows (see Theorem \ref{t5.4}):
\bt  \label{t1.3}  Let $M^{\ast}=\CP \setminus V$ be the complement of an essential hyperplane arrangement. Then, for any $k\leq n-1$,  \be   
\bigcup_{i=0}^{k}  \V_{i}(M)= \bigcup_{i=0}^{k}  \V_{i}(M^{\ast})  \subset  \bigcup_{S\subset V_{1}} \lbrace\V(M^{\ast}_{x_{S}}) \times V(\prod_{i \notin I_{x_S}} t_{i}=1)\rbrace,
\ee where the union is over the strata of $V_{1}$ with $\dim_{\C} S\geq n-k-1$, and $x_{S}\in S$ is arbitrary.  
\et

In the course of proving Theorem \ref{t1.2}, we also investigate (and obtain divisibility results for) the Alexander and characteristic varieties of the space $$\U:=\CP \setminus (V\cup H),$$
where $H$ is any hyperplane in $\CP$, which is in general position with respect to the hypersurface $V$. These are then related to $\V_{i}(M)$ and $\V_{i}(M^{\ast})$ via the $n$-homotopy equivalence $\U \hookrightarrow M$, see Proposition \ref{p5.2}.
 
 \medskip

As an application of Theorem \ref{t1.2}, we obtain the following characterization of the topological Euler characteristic of $F$ and $M^{\ast}$ (see Corollary \ref{c5.4}), which generalizes several known results obtained by different techniques:
\bc  Assume $r\geq 2$ and $V$ satisfies one of the following assumptions:
\begin{enumerate}
\item[(a)]  $V$ is an essential hypersurface, i.e.,  $\bigcap_{i=1}^{r} V_{i}= \emptyset$.
\item[(b)]  One of the irreducible components of $V$, say $V_{1}$, is smooth and transversal to $\bigcap_{i=2}^{r} V_{i}$.
\item[({c})] $V$ is a normal crossing divisor along any point of $\bigcap_{i=1}^{r} V_{i}$.
\end{enumerate} Then 
\be
  (-1)^{n}\chi(F)= (-1)^{n} d \cdot \chi(M^{\ast})\geq 0.
   \ee
\ec

\subsection{Summary and applications}

The paper is organized as follows.

In Section \ref{S2}, we investigate the relationship between the homology and resp. cohomology Alexander varieties. In Section \ref{S3}, we introduce the Sabbah specialization complex associated to an affine hypersurface, and compute 
 the stalks of its cohomology sheaves in terms of the corresponding local cohomology Alexander modules.
The proofs of the above-mentioned divisibility results of Theorem \ref{t1.2} are carried out in Sections \ref{S4} and \ref{S5}. In particular, we derive vanishing results for the homology Alexander polynomials in Section \ref{s5.2}. In Section \ref{s5.3}, we focus on the special case of hyperplane arrangements.

Sections \ref{S6} and \ref{S7} are devoted to applications of these divisibility results. 

In section \ref{S6}, we first recall the definition of the Dwyer-Fried invariants. Dwyer and Fried \cite{DF} observed that the homological finiteness properties of free abelian covers are completely determined by the Alexander varieties.   Moreover, it was noted in \cite{PS10} that in the previous sentence one can replace the Alexander varieties by the characteristic varieties. By combining this fact with our divisibility results, we recast several finiteness  and divisibility results for classical (infinite cyclic) Alexander modules of hypersurface complements, e.g., Libgober's results for  hypersurfaces with only isolated singularities, including at infinity \cite{L94}, or the second-named author's results for  hypersurfaces in general position at infinity \cite{Max}. 
We also obtain a divisibility result for the classical Alexander invariants of essential hyperplane arrangements, which (to our knowledge) seems to be new (see Remark \ref{r6.1}):
\bp
Assume that $V$ is an essential hyperplane arrangement in $\CP$ with $r$ irreducible components. Choose one of the irreducible components of $V$, say $V_{1}$, as the hyperplane at infinity. Then the classical (infinite cyclic) Alexander modules of $M^{\ast}$ are  torsion up to degree $n-1$. If $\omega_{i}(t)$  ($i\leq n-1$) denotes the corresponding $i$-th Alexander polynomial of $M^{\ast}$, then the roots of $\omega_{i}(t)$ have order $(r-1)!$.
\ep

In Section  \ref{S7}, we first introduce the resonance varieties $ \re^{i}(-)$ and recall the tangent cone inclusion (cf. \cite{L02}), which establishes a nice relationship between the characteristic varieties and the resonance varieties. 
It is known that for the special case of hyperplane arrangements, the tangent cone inclusion becomes an equality. Using this fact, we translate our divisibility results for the Alexander and characteristic varieties of hyperplane arrangement complements into similar results for the corresponding resonance varieties (see Theorem \ref{t7.3}): 
\bt   Let $M^{\ast}$ be the complement of an essential hyperplane arrangements. Then, for any $k\leq n-1$, we have:
\be 
 \bigcup_{i=0}^{k}  \re^{i}(M)= \bigcup_{i=0}^{k}  \re^{i}(M^{\ast}) \subset  \bigcup_{S\subset V_{1}} \lbrace \re(M^{\ast}_{x_{S}}) \times V(\sum _{i\notin I_{x_S}} z_{i}=0) \rbrace ,
\ee
where the union is over the strata of $V_{1}$ with $\dim_{\C} S\geq n-k-1$ and $x_{S}\in S$ is arbitrary.
\et
\noindent Here, $\re^{i}(M)$ and $\re^{i}(M^{\ast})$ denote the $i$-th resonance varieties of $M$ and $M^{\ast}$,  respectively, while for a finite connected $n$-dimensional  CW-complex $X$, we set (see Definition \ref{d7.0}): $$\re(X)=\bigcup_{i=1}^n \re^i(X).$$ Note that the resonance varieties $\re^{i}(M)$ and $\re^{i}(M^{\ast})$ are both homogeneous sub-varieties of the affine space $\C^{r}$, with coordinates $z=(z_{1},\cdots,z_{r})$.

\medskip

\noindent{\it Convention:} Unless otherwise specified, all homology and cohomology groups will be assumed to have $\C$-coefficients.

\medskip

\noindent{\it Disclaimer:} We assume reader's familiarity with the standard homological algebra language, as well as with basic knowledge of the theory of constructible sheaf complexes (and in particular, of perverse sheaves).

\medskip

\textbf{Acknowledgments.} 
We thank Nero Budur for valuable comments on an earlier version of this paper.
The first named author is supported by China Scholarship Council (file No.
201206340046). He thanks the Mathematics Department at the University of Wisconsin-Madison for hospitality during the preparation of this work. The second-named author was partially supported by grants from NSF (DMS-1304999), NSA (H98230-14-1-0130), and by a grant of the Ministry of National
Education (CNCS-UEFISCDI project number PN-II-ID-PCE-2012-4-0156).


\section{Alexander varieties and Characteristic varieties}\label{S2}
In this section, we introduce three types of jumping loci associated to a space and
establish a relation between homology Alexander varieties and cohomology Alexander varieties.

\subsection{Algebraic preliminaries}\label{s2.1}
Let $\Gamma_{r}=\C[t_{1},t_{1}^{-1},\cdots,t_{r},t_{r}^{-1}]$ be the Laurent polynomial ring in $r$ variables. $\Gamma_{r}$ is a regular Noetherian domain, and in particular it is factorial. If $A$ is a finitely generated $\Gamma_{r}$-module, the {\it support locus of $A$} is defined by the annihilator ideal of $A$, i.e.,  
\be \Su(A)=V(ann(A)).\ee 
By construction, $\Su(A)$ is a Zariski closed subset of the $r$-dimensional affine torus $(\C^{\ast})^{r}$, where $(\C^{\ast})^{r}$ is the scheme defined by $\Gamma_{r}$. The support locus $\Su(A)$ is called {\it proper}, if it has at least codimension $1$ in $(\C^{\ast})^{r}$.

 For $\lambda=(\lambda_{1},\cdots, \lambda_{r})\in (\C^{\ast})^{r}$, let $m_{\lambda}$ be the corresponding maximal ideal in $\Gamma_{r}$. For a $\Gamma_{r}$-module $A$, we denote by $A_{\lambda}$ the localization of $A$ at the maximal ideal $m_{\lambda}$. For $A=\Gamma_{r}$, we use the simpler notation $\Gamma_{\lambda}$. Then, if $A$ is of finite type, we have that \be \label{2.1}
 \Su(A)=\lbrace \lambda \in (\C^{\ast})^{r}\mid A_{\lambda}\neq 0 \rbrace .
 \ee 

Let $P$ be a prime ideal in $\Gamma_{r}$, of height $1$. Since $\Gamma_{r}$ is a regular Noetherian domain, the localization of $\Gamma_{r}$ at $P$, denoted by $\Gamma_{P}$, is a regular local ring of dimension $1$. So $\Gamma_{P}$ is a discrete valuation ring, hence a principal ideal domain. A prime ideal $P$ of height $1$ in $\Gamma_{r}$ is principal, and we let $\Delta({P})$ denote the generator of $P$, which is well-defined up to multiplication by units of $\Gamma_{r}$.

Assume that $\Su (A)$ is proper in $(\C^{\ast})^{r}$, i.e., $ann(A)\neq 0$.   Then, if $V(P) \nsubseteq \Su(A)$, $A_{P}=0$. On the other hand, if $V(P)\subset \Su(A)$, then $ann(A) \subset P$ and  $ann(A_{P})=\lbrace ann(A)\rbrace_{P} \neq 0$, hence $A_{P}$ is a torsion $\Gamma_{P}$-module. This shows that $A_{P}$ has finite length as a $\Gamma_{P}$-module, which is denoted by $lg(A_{P})$. The {\it characteristic polynomial of $A$} is defined by \be \Delta(A)=\prod_{P} \Delta({P})^{lg(A_{P})},\ee 
where the product is over all the prime ideals $P$ in $\Gamma_{r}$ of height $1$ such that $V(P)\subset \Su(A)$. Since $\Su A$ is proper, this product is indeed a finite product. The prime factors of $\Delta(A)$ are in one-to-one corespondence with the codimension one irreducible hypersurfaces of $(\C^{\ast})^{r}$ contained in $\Su(A)$. Note that if ${codim} \ \Su(A)>1$ then $\Delta(A)=1$. If $\Su(A)=(\C^{\ast})^{r}$, then $\Delta(A)=1$ by convention. 

\br Sabbah used the alternating product of the above characteristic polynomials to define the Zeta function of a $\Gamma_{r}$-module complex, see \cite[Definition 2.1.12]{Sab}.
The support of the $\Gamma_{r}$-module $A$ can also be defined by the first Fitting ideal of $A$, see \cite{DM}. 
\er

\bl  \label{l2.1} If $A \to B \to C$ is an exact sequence of finitely generated $\Gamma_{r}$-modules, then  $$\Su(B) \subset \Su(A) \cup \Su({C}).$$ Moreover, if $\Su(A)$ and $\Su({C})$ are proper in $(\C^{\ast})^{r}$, then so is $\Su(B)$, and $\Delta(B)$ divides $\Delta(A)\times \Delta(C)$.
 \el
\begin{proof}  Using the exactness of localization (see \cite{Wei}, p.76), for any $\lambda \in (\C^{\ast})^{r}$ as above we get an exact sequence $A_{\lambda} \to B_{\lambda}  \to C_{\lambda} $. If $\lambda \notin \Su(A) \cup \Su(C)$, then $A_{\lambda}=0=C_{\lambda}$. Thus $B_{\lambda}=0$, i.e., $\lambda\notin \Su(B)$.

 Let $P$ be a prime ideal in $\Gamma_{r}$ of height $1$. By the exactness of localization, we get an exact sequence $A_{P} \to B_{P}  \to C_{P} $.   
Assuming that $\Su(A)$ and $\Su({C})$ are proper in $(\C^{\ast})^{r}$, it follows that $lg(B_{P})\leq lg(A_{P})+ lg(C_{P}) $. So $\Delta(B)$ divides $\Delta(A)\times \Delta(C)$.
\end{proof}
\bl   \label{l2.2}  \cite[Proposition 4(c)]{Ser} If $A$ and $B$ are finitely generated $\Gamma_{r}$-modules, then 
 $$\Su(A\otimes_{\Gamma_{r}} B)=\Su(A)\cap \Su(B).$$
\el


\subsection{Alexander varieties}\label{s2.2}
Assume that $X$ is a finite connected $n$-dimensional CW complex with $\pi_{1}(X)=G$. Let $\nu: G\twoheadrightarrow \Z^{r}$ be an epimorphism, and consider the corresponding  free abelian cover $X^{\nu}$ of $X$. The group of covering transformations of $X^{\nu}$ is isomorphic to $\Z^{r}$ and acts on the covering space. By choosing fixed lifts of the cells of $X$ to $X^{\nu}$, we obtain a free basis for the chain complex  $C_{\ast}(X^{\nu},\C)$ of  $\Gamma_{r}$-modules.
 
\br \label{r2.6} (\cite[Remark 2.2]{DM})   Though the ring $\Gamma_{r}$ is commutative, it should be regarded as a quotient
ring of $\C[\pi_{1}(X)]$, which is non-commutative in general. Hence, one should
be careful to distinguish the right from the left $\Gamma_{r}$-modules. The ring $\Gamma_{r}$ has a natural involution denoted by an overbar, sending each $t_{i}$ to  $\overline{ t_{i}}:= t_{i}^{-1}$. To a left $\Gamma_{r}$-module $A$, we associate the right $\Gamma_{r}$-module $\overline{A}$, with the same underlying abelian group but with the $\Gamma_{r}$-multiplication given by $a\cdot c\longmapsto \overline{c}\cdot a$, for $a\in A$ and $c\in \Gamma_{r}$. 
In this paper, we regard $C_{\ast}(X^{\nu},\C)$  as a complex of right $\Gamma_{r}$-modules.
Conversely, to a right $\Gamma_{r}$ module $A$, we can associate a corresponding left module $\overline{A}$ by using the involution.    Moreover, it can be seen that $A= \overline{(\overline{A})}$.
 \er
 
 Let us now define $$\overline{\Su(A)}:=\lbrace (t_{1},\cdots, t_{r})\in (\C^{\ast})^{r} \mid (t_{1}^{-1},\cdots, t_{r}^{-1}) \in \Su(A)\rbrace$$ and $$\overline{\Delta(A)}(t_{1},\cdots,t_{r}):=\Delta(A)(t_{1}^{-1},\cdots,t_{r}^{-1}).$$ Then $\Su(\overline{A})=\overline{\Su(A)}$ and $\Delta(\overline{A})=\overline{\Delta(A)}$.
 
\bd \label{d2.1}  The {\it $i$-th homology Alexander module} $A_{i}(X,\nu)$ of $X$ associated to the epimorphism $\nu$ is by definition the $\Gamma_{r}$-module $ H_{i}(X^{\nu},\C)$. Similarly, the {\it $i$-th cohomology Alexander module} $A^{i}(X,\nu)$ of $X$ is by definition the $i$-th cohomology module of the dual complex $\Hom_{\Gamma_{r}} ({C_{\ast}(X^{\nu},\C)}, \Gamma_{r})$, where $\Gamma_{r}$ is considered here with  the induced right $\Gamma_{r}$-module structure, as in Remark \ref{r2.6}. 
\ed

Since $X$ is a finite $n$-dimensional CW complex, it is clear by definition that 
$$A_{i}(X,\nu)=0=A^{i}(X,\nu) \text{  for  } i>n.$$
 
 \bd \label{d2.2} The {\it $i$-th homology Alexander variety of $X$ associated to the epimorphism $\nu$} is the support locus of the annihilator ideal of  $A_{i}(X,\nu)$: 
 $$\W_{i}(X,\nu):=\Su(A_{i}(X,\nu)) ,$$ and, the corresponding {\it $i$-th homology Alexander polynomial} is $\Delta_{i}(X,\nu):=\Delta(A_{i}(X,\nu))$.
   Similarly, the {\it $i$-th cohomology Alexander variety of $X$ associated to the epimorphism $\nu$} is the support locus of the annihilator ideal of  $A^{i}(X,\nu)$: 
 $$\W^{i}(X,\nu):=\Su(A^{i}(X,\nu)) ,$$ and, the corresponding {\it $i$-th cohomology Alexander polynomial} is $\Delta^{i}(X,\nu):=\Delta(A^{i}(X,\nu))$.
  \ed
  
 Since $X$ is a finite CW-complex, the chain complex $C_{\ast}(X^{\nu},\C)$ is of finite type, so the corresponding (homology or cohomology) Alexander varieties are well-defined. The (homology or cohomology) Alexander varieties and (homology or cohomology) Alexander polynomials are homotopy invariants of $X$. 
  
\br For the one-variable case, the homology Alexander polynomials coincide with the classical Alexander polynomials, as studied e.g., in  \cite{L94,Max}, while the homology Alexander varieties can be identified with the sets of roots of the corresponding classical Alexander polynomials.\er

  If $\nu$ is the composition map: $\pi_{1}(X) \twoheadrightarrow H_{1}(X,\Z) \twoheadrightarrow \text {Free  }H_{1}(X,\Z)$, where the first map is the Hurewicz (or abelianization) map and the second map just kills the torsion part, then we simply write $A_{i}(X):=A_{i}(X,\nu)$, $\W_{i}(X):=\W_{i}(X,\nu)$, $\Delta_{i}(X):=\Delta_{i}(X,\nu)$, and similar notations are used for the cohomology versions. Moreover, if $H_{1}(X,\Z)$ is free, this corresponding cover is called the {\it universal abelian  cover}, and the covering map is denoted by $ab$. 
 
 Here we only give the definition of these Alexander invariants with $\C$-coefficients, but see \cite[Section 3.2]{PS10} for arbitrary field coefficients.  


\subsection{Homology versus cohomology Alexander varieties}\label{s2.3}
In this section, we employ here the spectral sequence approach from \cite[section 2.3]{DM} in order to relate the homology and resp. cohomology Alexander varieties. More precisely,  $A^{\ast}(X)$ can be computed from  $A_{\ast}(X)$ by using the Universal Coefficient spectral sequence: 
\be \label{ss2.1}
\Ext^{q}_{\Gamma_{r}}(A_{s}(X), \Gamma_{r}) \Longrightarrow A^{s+q}(X) .
\ee
 Using the exactness of localization, we then get the following localized spectral sequence for any $\lambda \in (\C^*)^r$: 
 \be  \label{ss2.2}
  \Ext^{q}_{\Gamma_{\lambda}}({A_{s}(X)_{\lambda}}, \Gamma_{\lambda}) \Longrightarrow A^{s+q}(X)_{\lambda} .
 \ee
 The next proposition is a direct consequence of \cite[Proposition 2.3]{DM}.
 
 \bp \label{p2.1} For $X$ a finite connected CW-complex and any $k\geq 0$, we have:  
 \be  \label{2.2}
 \bigcup _{i=0}^{k}  \W^{i}(X) \subset  \bigcup_{i=0}^{k} { \W_{i}(X)}.
 \ee
 \ep
 \begin{proof}
 If $\lambda \notin  \bigcup _{i=0} ^{ k} { \W_{i}(X)}$, then we have by (\ref{2.1}) that ${ A_{i}(X)_{\lambda}}=0$ for all  $0 \leq i\leq k$. Then the spectral sequence (\ref{ss2.2}) yields that $A^{i}(X)_{\lambda}=0$ for $i\leq k$. Hence $\lambda \notin \bigcup _{i=0} ^{ k} \W^{i}(X)$. 
 \end{proof}
 
Since $C_{\ast}(X^{ab},\C)$ is a complex of finitely generated free $\Gamma_{r}$-modules, we have the following $\Gamma_{r}$-module complex isomorphism: 
 $$\Hom_{\Gamma_{r}} ( \overline{\Hom_{\Gamma_{r}} ({ C_{\ast}(X^{ab},\C}), \Gamma_{r} )},\Gamma_{r}) \cong  \overline{C_{\ast}(X^{ab},\C)}.$$
  So, similarly, we also have a Universal Coefficient spectral sequence computing $A^{\ast}(X)$ from $A_{\ast}(X)$:
 \be \label{ss2.3}
 \Ext^{q}_{\Gamma_{s}}(\overline{ A^{n-s}(X)}, \Gamma_{r}) \Longrightarrow \overline{A_{n-s-q}(X)} 
\ee
 and a corresponding localized spectral sequence for any $\lambda \in (\C^*)^r$: 
 \be  \label{ss2.4}
 \Ext^{q}_{\Gamma_{\lambda}}(\overline{ A^{n-s}(X)_{\lambda}}, \Gamma_{\lambda}) \Longrightarrow \overline{A_{n-s-q}(X)}_{\lambda}.
 \ee

 \bp \label{p2.2} For $X$ a finite connected CW-complex of dimension $n$, and any $k\leq n$, we have: 
 \be  \label{2.3}
 \bigcup _{i=k}^{n} \W_{i}(X) \subset  \bigcup _{i=k}^{n} { \W^{i}(X)}.
 \ee
 Moreover, \be \label{2.4}
  \bigcup _{i=0}^{n} \W_{i}(X) =  \bigcup _{i=0}^{n}  {\W^{i}(X)}.
 \ee
 \ep
 \begin{proof}
 The first claim follows from the spectral sequence (\ref{ss2.4}) by an argument similar to that of Proposition \ref{p2.1}. The second claim follows from (\ref{2.2}) and (\ref{2.3}).
 \end{proof}
 
\br \label{r2.1} Since $A_{0}(X)=\Gamma_{r}/m_{\fst}$, where $\fst=(1,\cdots,1)$, it follows from  \cite[Proposition 2.3]{DM}  that $A^{0}(X)=0$.
\er 

\br \label{r2.2}  Since $X$ is a finite $n$-dimensional CW-complex, $A_{n}(X)$ is torsion-free as a $\Gamma_{r}$-module.  So, if $A_{n}(X)\neq 0$, then $\W_{n}(X)=(\C^{\ast})^{r}$.  Moreover, (\ref{2.3}) yields that 
$  \W_{n}(X) \subset { \W^{n}(X)} $, hence $\W^{n}(X)=(\C^{\ast})^{r}$. In particular, $\Delta_{n}(X)=\Delta^{n}(X)=1$. 
\er 
 
Next, we investigate the relationship between the Alexander polynomials $\Delta_{\ast}(X)$ and $\Delta^{\ast}(X)$.

\bp \label{p2.3}  Assume that $\bigcup_{i=0}^{k} \W_{i}(X)$ is proper in $(\C^*)^r$. Then, for any $i\leq k$, \be \label{2.5}
{\Delta_{i-1}(X)}= \Delta^{i}(X).
\ee
Therefore, there is an exact correspondence between the codimension 1 irreducible hypersurfaces of $(\C^{\ast})^{r}$ contained in ${\W_{i-1}(X)}$ and $\W^{i}(X)$, respectively, for $i\leq k$.
\ep
\begin{proof}
Let $P$ be a prime ideal in $\Gamma_{r}$ of height $1$. By the exactness of localization at $P$, we get from (\ref{ss2.1}) a localized spectral sequence: 
\be   \label{ss2.5}
E_{2}^{s,q}=  \Ext^{q}_{\Gamma_{P}}({ A_{s}(X)_{P}}, \Gamma_{P}) \Longrightarrow A^{s+q}(X)_{P} .
 \ee
 $\Gamma_{P}$ is a discrete valuation ring, hence a principal ideal domain. So $E_{2}^{s,q}=0$ for $q>1$. Note also that $E_{2}^{s,0}=\Ext^{0}_{\Gamma_{P}}({ A_{s}(X)_{P}}, \Gamma_{P})$ is always free, and $E_{2}^{s,1}=\Ext^{1}_{\Gamma_{P}}({ A_{s}(X)_{P}}, \Gamma_{P})$ is always torsion as a $\Gamma_{P}$-module. So the spectral sequence (\ref{ss2.5}) degenerates at $E_{2}$. 
Since $\bigcup_{i=0}^{k} \W_{i}(X)$ is proper, we have that 
$$\Ext^{0}_{\Gamma_{r}}({ A_{s}(X)},\Gamma_{r})=\Hom_{\Gamma_{r}}({ A_{s}(X)},\Gamma_{r})=0, \text{ for any } s\leq k.$$
Then the spectral sequence (\ref{ss2.5}) yields that $\Ext^{1}_{\Gamma_{P}}({ A_{i-1}(X)_{P}}, \Gamma_{P}) \cong A^{i}(X)_{P}$, for $i\leq k$. The claim follows now from the Universal Coefficient Theorem for the principal ideal domain $\Gamma_{P}$.
\end{proof}

\br  The results stated in this subsection hold more generally, for the Alexander varieties and Alexander polynomials associated to any epimorphism $\nu: G\twoheadrightarrow \Z^{r}$. 
\er 


\subsection{Characteristic varieties}\label{s2.4}
Let $X$ be, as before, a finite $n$-dimensional connected  CW-complex with $\pi_{1}(X)=G$. Let $\C^{\ast}$ be the multiplicative group of non-zero complex numbers. Then the group of $\C$-valued  characters, $ \Hom(G,\C^{\ast})$, is a commutative, affine algebraic group. Each character $\rho \in \Hom(G,\C^{\ast})$ defines a rank one local system on $X$, denoted by $\K_{\rho}$.
 \bd \label{d2.3} The {\it $i$-th characteristic variety of $X$} is the Zariski closed subset of the character group $ \Hom(G,\C^{\ast})$ defined as: 
$$\V_{i}(X)=\lbrace \rho\in \Hom(G,\C^{\ast}) \mid \dim_{\C} H_{i}(X,\K_{\rho})\neq 0 \rbrace.$$ 
\ed
\br
$\V_{i}(X)$ is filtered by the closed subsets
$$\V^k_{i}(X)=\lbrace \rho\in \Hom(G,\C^{\ast}) \mid \dim_{\C} H_{i}(X,\K_{\rho})\geq k \rbrace,$$ for all $k>0$.
\er
\br \label{r2.4} Some references define the characteristic varieties by using cohomology rather than homology:
 $$\V^{i}(X)=\lbrace \rho\in \Hom(G,\C^{\ast}) \mid \dim_{\C} H^{i}(X,\K_{\rho})\neq 0 \rbrace.$$ 
Let $\overline{\rho}$ denote the inverse of $\rho$ in  the character group $\Hom(G,\C^{\ast})$, then \be \label{2.8}
 H^{i}(X,\K_{\overline{\rho}}) \cong \Hom_{\C}(H_{i}(X,\K_{\rho}),\C).
 \ee  Therefore, $\V_{i}(X)=\overline{\V^{i}(X)}$, where $\overline{\V^{i}(X)}=\lbrace \rho\in  \Hom(G,\C^{\ast}) \mid \overline{\rho} \in \V^{i}(X)\rbrace$.
\er

The characteristic varieties $\V_{i}(X)$ are homotopy invariants of $X$, see \cite[Lemma 2.4]{Su14}.
\medskip 
 
If $\nu: \pi_1(X)=G\twoheadrightarrow \Z^{r}$ is an epimorphism, then $\nu$ induces an embedding
$$\nu^{\ast}: (\C^{\ast})^{r}=\Hom(\Z^{r},\C^{\ast}) \hookrightarrow \Hom(G,\C^{\ast}).$$ 
The following relationship between the characteristic varieties and the homology Alexander varieties was established in \cite{PS10}:
\bt \label{t2.1}   (\cite[Theorem 3.6]{PS10}) For any $k\geq 0$ and any epimorphism $\nu: G\twoheadrightarrow \Z^{r}$,
\be \label{2.5}
\bigcup_{i=0}^{k} \W_{i}(X,\nu)   =im(\nu^{\ast}) \cap \left(\bigcup_{i=0}^{k} \V_{i}(X)\right) .
\ee
\et 

Let $\Hom(G,\C^{\ast})^{0}$  denote the identity component of the algebraic group $\Hom(G,\C^{\ast})$, with $G=\pi_1(X)$ as before. Set:
$$\V(X)=\bigcup _{i=0}^{n} \V_{i}(X) \ \text { and } \ \W(X)=\bigcup _{i=0}^{n} \W_{i}(X).$$
Then the following result follows at once from Theorem \ref{t2.1}:
\bp \label{p2.4} If $X$ is a finite connected CW-complex of dimension $n$, then: 
\be \label{2.7}
 \W(X)=\V(X) \cap \Hom(G,\C^{\ast})^{0}.
\ee 
 \ep

\br \label{r2.3} It is known that for any character $\rho$ of $\pi_1(X)$, the following identity holds: $\chi(X,\K_{\rho})=\chi(X)$ (see \cite[Proposition 2.5.4]{D2}). Hence, if $\chi(X)\neq 0$, we have that $\V(X)=\Hom(G,\C^{\ast})$ and $\W(X)=\Hom(G,\C^{\ast})^{0}$.  So, $\V(X)$ and $\W(X)$ are interesting to study only when $\chi(X)=0$.  
\er


\section{Sabbah Specialization complex and local Alexander modules}\label{S3}
In this section, we recall the definition of the {\it Sabbah specialization complex} associated to an affine hypersurface, and compute the stalks of its cohomology sheaves in terms of  the local cohomology Alexander modules. 

\medskip

For any complex algebraic variety $X$ and any commutative ring $R$, we denote by  $D^{b}_{c}(X,R)$ the derived category of bounded cohomologically $R$-constructible complexes of sheaves on $X$. For a quick introduction to derived categories, the reader is advised to consult \cite{D2}.

\bd  For $\F \in D_{c}^{b}(X,\Gamma_{r})$ and a point $x\in X$, the {\it $i$-th support of $\F$ at $x$} is defined by  
$$\Su^{i}_x (\F) := \Su (\mathcal{H} ^{i}(\F)_{x}) \subset (\C^{\ast})^{r},$$ 
and the {\it $i$-th characteristic polynomial of $\F$ at $x$} is defined by
$$\delta^{i}(\F_{x}) :=\Delta (\mathcal{H} ^{i}(\F)_{x}).$$ 
The {\it support of $\F$ at $x$} is then set to be the union
$$\Su_{x} (\F) := \bigcup_{i} \Su^{i}_{x} (\F).$$ 
The {\it multi-variable monodromy zeta-function of $\F$ at $x$} is defined as:
$$Z(\F_{x}):=\prod_{i} \delta^{i}(\F_{x})^{(-1)^{i}}.$$
\ed


\subsection{Alexander modules of hypersurface complements}\label{s3.1}
Let $V$ be a hypersurface in $\CP$ defined by a degree $d$ reduced homogeneous polynomial $f$. Assume that $f=f_{1}\cdot \cdots \cdot f_{r}$, where $f_{i}$ are the irreducible factors of $f $, so $V_{i}=\lbrace f_{i}=0 \rbrace$ are the irreducible components of $V$. 

Choose a hyperplane $H \nsubseteq V$ in $\CP$ as the hyperplane at infinity. Without loss of generality, assume that $\CP$ has coordinates $[x_{0}:\cdots :x_{n}]$ and $H=\lbrace x_{0}=0\rbrace$. Then $\CN=\CP\setminus H$ has coordinates $(x_{1},\cdots,x_{n})$. Set $g_{i}=f_{i}(1,x_{1},\cdots,x_{n})$ and $g=\prod_{i=1}^{r}g_{i}$, so $\g=(g_{1},\cdots,g_{r})$ defines a polynomial map from $\CN$ to $\C^{r}$. Set \begin{center}
   $D_{i}=\lbrace g_{i}=0 \rbrace$,    $D=\bigcup^{r}_{i=1} D_{i} $,
 \end{center}
 and   $$\U= \CN\setminus D=\CP\setminus (V \cup H).$$
It is known that $H_{1}(\U,\Z)={\Z}^{r}$ is torsion free (\cite{D1},(4.1.3),(4.1.4)), generated by meridian loops $ \gamma_{i} $ about each irreducible component $ V_{i}$, $1\leq i\leq r $. If $ \gamma_{\infty} $ denotes the meridian about the hyperplane at infinity, then there is a relation in $H_{1}(\U,\Z)$ : $\sum_{i=1}^{r} d_{i}\gamma_{i}+\gamma_{\infty}=0 ,$ where $ d_{i}=\deg f_{i} $.

We can associate to the hypersurface complement $\U$ the homology (resp. cohomolgy) Alexander varieties $\W_{i}(\U)$ (resp. $\W^{i}(\U)$) and homology (resp. cohomolgy) Alexander polynomials $\Delta_{i}(\U)$ (resp. $\Delta^{i}(\U)$) as defined in Section 2.2. Note that $\U$ has the homotopy type of a finite $n$-dimensional CW complex. So $A_{i}(\U)=0=A^{i}(\U)$ for $i>n$, and $A_{n}(\U)$ is torsion-free as a $\Gamma_{r}$-module.
  
Define a local system $\K$ on $\U$ with stalk $\Gamma_{r}$ and representation of the fundamental group determined by the composition:
$$\pi_{1}(\U) \to H_{1}(\U,\Z) \to Aut(\Gamma_{r}),$$
where the first map is the abelianization, and the last map is defined by $\gamma_{i}\to t_{i}$, i.e., the action of $\gamma_i$ on $\Gamma_{r}$ is given by the multiplication with $t_{i}$. Let $\K^{\vee}$ be the dual local system, whose stalk at a point $x\in \U$ is $\K^{\vee}_{x}:=\Hom_{\Gamma_{r}}(\K_{x}, \Gamma_{r})$.

For $\F\in D^{b}_{c}(X,\Gamma_{r})$, let $\D \F$ denote its Verdier dual and $\overline{\F}$ denote its conjugation (i.e., $\overline{\F}$ is defined by composing all module structures with the involution). Then there is an isomorphism of local systems: $\K^{\vee} \simeq \overline{\K}$, see \cite[section 2.2]{DM}.

The Alexander modules of $\U$ are related to the above local system by the following $\Gamma_{r}$-module isomorphisms:
\begin{equation}\label{3.0}
H_{i}(\U,\K)\cong A_{i}(\U) \text{    and         } H^{i}(\U,\overline{\K})\cong A^{i}(\U)
\end{equation} 
 for all $i$.
Note that since $\U$ is smooth of complex dimension $n$, we also have that \begin{center}
$\mathcal{D} \K  = \overline{\K}[2n]$ \ and \ $\D (\D \K)  =\K$.
\end{center}

\br  \label{r3.1} The following isomorphism follows by an argument similar to the one used in the proof of \cite[Corollary 3.4]{Max}:\be \label{3.5}
A_{i}(\U)\cong H^{2n-i}_{c}(\U, \K).
\ee
The second named author proved this isomorphism in the one-variable case by using intersection homology, see \cite[Corollary 3.4]{Max}. Since the intersection homology theory is defined for any Noetherian commutative ring of finite
cohomological dimension (\cite[pp.68]{B}), this isomorphism also holds for the multi-variable Laurent polynomial ring  $\Gamma_{r}$.
  \er 
  
  
 \subsection{Sabbah specialization complex}\label{s3.2}
 The Sabbah specialization complex (see \cite{Sab}, and its reformulation in \cite{Bu}) can be regarded as a generalization of Deligne's nearby cycle complex. 

In our notations, consider the following commutative diagram of spaces and maps:\begin{center}
$\xymatrix{
D \ar[r]^{i}   & \CN  \ar[d]^{\g}   & \ar[l]_{l} \U \ar[d]^{\g}  & \ar[l]_{\pi} \U^{ab} \ar[d]^{\widehat{\g}} \\
                       & \C^{r}    & \ar[l] (\C^{\ast})^{r}                                           
            & \ar[l]_{\widehat{\pi}} \widehat{(\C^{\ast})^{r}} 
}$
\end{center}
where $\widehat{\pi}$ is the universal covering of $(\C^{\ast})^{r}$, and the right-hand square of the diagram is cartesian.

\bd (\cite[Definition 3.2]{Bu}) The Sabbah specialization complex functor  of $\g$ is defined by \begin{center}
$\psi_{\g} =i^{\ast}Rl_{\ast} R\pi_{!} (l \circ \pi)^{\ast} : D^{b}_{c}(\CN,\C) \rightarrow D^{b}_{c}(D,\Gamma_{r})$, 
\end{center} and we call $\psi_{\g}\C_{\CN}$ the {\it Sabbah specialization complex}. 
\ed
 
\begin{remark}  (\cite[Remark 3.3]{Bu}) The Sabbah specialization complex can be viewed as a generalization of Deligne's nearby cycles. In fact, when $r = 1$, $\psi_{\g}\F$ as defined here equals $\psi_{\g}\F[-1]$ as defined by Deligne (see \cite{Br}, page 13), where $R\pi_{!}$ is replaced by $R\pi_{\ast}$. 
\end{remark} 
 For short, in the following we write $\C$ for the constant sheaf $\C_{\CN}$ on $\CN$.
 
\begin{lem} \label{l3.1} \cite[Lemma 3.4]{Bu}
   We have the following distinguished triangle in $D^{b}_{c}(\CN,\Gamma_{r})$: \be\label{3.1}
l_{!}\K \rightarrow Rl_{\ast}\K \rightarrow   i_{!}\psi_{\g}\C \overset{[1]}{\rightarrow}.
\ee In particular,  $\psi_{\g}\C= i^{\ast}Rl_{\ast}\K$.
\end{lem}


\subsection{Local Alexander modules}\label{s3.3}
Let $x$ be a point in $D$ and $B_{x}$ be a small open ball centered at $x$ in $\CN$. Set $\U_{x}=\U\cap B_{x}$. Denote by $I_{x}=\lbrace i \mid g_{i}(x)= 0 \rbrace$ and $r_{x}=\# \vert I_{x} \vert$. It is known that $H_{1}(\U_{x},\Z)$ is torsion free. The map $H_{1}(\U_{x},\Z)\to H_{1}(\U,\Z)$ induced by the natural inclusion can be viewed as a composition of the following two maps: \be \label{3.2}
H_{1}(\U_{x},\Z)\overset{}{\twoheadrightarrow} \Z^{r_{x}}  \hookrightarrow  \Z^{r}= H_{1}(\U,\Z)
\ee
where $\Z^{r_{x}}$  is the free abelian group generated by the meridian loops $\lbrace \gamma_{i} \rbrace_{i\in I_{x}}$. The first map in (\ref{3.2}) is surjective and the second map is  injective. (In the special case of hyperplane arrangements, the first map in (\ref{3.2}) is always an isomorphism.)

 Let $\nu_{x}$ denote the canonical epimorphism: 
 \be
 \pi_{1}(\U_{x}) \overset{\nu_{x}}{\twoheadrightarrow} \Z^{r_{x}}
\ee 
defined by the composition 
$$ \pi_{1}(\U_{x}) \twoheadrightarrow H_{1}(\U_{x}\Z) \twoheadrightarrow \Z^{r_{x}}$$ of the first map in (\ref{3.2}) with the abelianization map. 
Let $\V(\U_{x})$ denote as before the union of the characteristic varieties of $\U_{x}$. Then Theorem \ref{t2.1} yields that $\W(\U_{x},\nu_{x})= \V(\U_{x})\cap im(\nu_{x}^{\ast})$. The embedding $\nu_{x}^{\ast}$ identifies the variables along branches that come from the same global irreducible component. For a concrete situation, see Examples \ref{ex3.4} and \ref{ex3.5}.

 Let $\Gamma_{r_{x}}$ be the sub-ring of $\Gamma_{r}$ with variables $\lbrace t_{i} \rbrace_{i\in I_{x}}$, and $\Gamma_{r-r_{x}}$  be the sub-ring of $\Gamma_{r}$ in the remaining variables. Then $\Gamma_{r}=\Gamma_{r_{x}}\otimes_{\C}\Gamma_{r-r_{x}}$.
\bd \label{d3.0} The {\it uniform Alexander varieties} of $\U_{x}$ are defined by 
\begin{center} $\W_{i}^{unif}(\U_{x},\nu_{x})=\W_{i}(\U_{x},\nu_{x})\times (\C^{\ast})^{r-r_{x}}$ \ and \ $\W^{unif}(\U_{x},\nu_{x})=\W(\U_{x},\nu_{x})\times (\C^{\ast})^{r-r_{x}}$,\end{center} where $(\C^{\ast})^{r-r_{x}}$ is the scheme corresponding to $\Gamma_{r-r_{x}}$.
\ed

 The characteristic varieties of the local complement $\U_{x}$ have been studied by Libgober in \cite{L09} and Budur-Wang in \cite{BW15}. In particular, $\V(\U_{x})$ is a finite union of torsion-translated subtori. By combining this property with the fact that the Alexander variety $\W(\U_{x},\nu_{x})= \V(\U_{x})\cap im(\nu_{x}^{\ast})$ can be defined over $\Q$, we get that the involution does not change $\W(\U_{x},\nu_{x})$. Moreover, if $\W(\U_{x},\nu_{x})$ is proper, then the involution keeps $\Delta^{i}(\U_{x},\nu_{x})=\Delta_{i-1}(\U_{x},\nu_{x})$ unchanged.

 \bp \label{p3.1}  With the above assumptions and notations, for any $x\in D$, we have:
 $$\Su^{i}_{x} (\psi_{\g}\C) ={ \W^{i}(\U_{x},\nu_{x})}\times (\C^{\ast})^{r-r_{x}};$$ $$\Su_{x}(\psi_{\g}\C)=\W^{unif}(\U_{x},\nu_{x}).$$ Moreover, $\Su_{x}(\psi_{\g}\C)$ is proper in $(\C^{\ast})^{r}$, so, for any $i\geq 0$,
 \be \label{3.3}
\delta^{i}(\psi_{\g}\C)_{x}= \Delta^{i}(\U_{x},\nu_{x})=\Delta_{i-1}(\U_{x},\nu_{x}). 
\ee
 \ep
\begin{proof} By Lemma \ref{l3.1}, there exist $\Gamma_{r}$-module isomorphsims: $\mathcal{H} ^{i}(\psi_{\g}\C)_{x}\cong H^{i}(\U_{x},\K\vert_{\U_{x}})$, for all $i$.  On the other hand, we have: 
\be  \label{3.4}
C^{\ast}(\U_{x},\overline{ \K}\vert_{\U_{x}})\cong \Hom_{\Gamma_{r_{x}}}(C_{\ast}(\U_{x}^{\nu_{x}}),\Gamma_{r})\cong \Hom_{\Gamma_{r_{x}}}(C_{\ast}(\U_{x}^{\nu_{x}}),\Gamma_{r_{x}})\otimes_{\C} \Gamma_{r-r_{x}},
\ee 
where the first isomorphism follows from the definition of the cohomology of a local system.
In fact, since $C_{\ast}(\U_{x}^{\nu_{x}})$ is isomorphic to a direct sum of copies of the ring $\Gamma_{r_{x}}$, we only need to show that (\ref{3.4}) holds for $C_{\ast}(\U_{x}^{\nu_{x}})= \Gamma_{r_{x}}$. As already stated, the first isomorphism in (\ref{3.4}) follows from the definition of the local system cohomology. Moreover, 
$$\Hom_{\Gamma_{r_{x}}}(\Gamma_{r_{x}},\Gamma_{r})\cong \Gamma_{r} \cong \Gamma_{r_{x}}\otimes_{\C} \Gamma_{r-r_{x}} \cong  \Hom_{\Gamma_{r_{x}}}(\Gamma_{r_{x}},\Gamma_{r_{x}})\otimes_{\C} \Gamma_{r-r_{x}},$$ where all the isomorphisms are $\Gamma_r$-module isomorphisms. Hence (\ref{3.4}) follows.

So $\mathcal{H} ^{i}(\psi_{\g}\C)_{x}\cong \overline{ A^{i}(\U_{x},\nu_{x})}\otimes_{\C} \Gamma_{r-r_{x}} $, and $ann(\mathcal{H} ^{i}(\psi_{\g}\C)_{x})= \langle ann (\overline{ A^{i}(\U_{x},\nu_{x})})\rangle,$
 where $\langle ann(\overline{ A^{i}(\U_{x},\nu_{x})})\rangle$ is the ideal in $\Gamma_{r}$ generated by $ann(\overline{ A^{i}(\U_{x}, \nu_{x})})$.   Therefore, $\Su^{i}_{x} (\psi_{\g}\C) =\overline{\W^{i}(\U_{x},\nu_{x})}\times (\C^{\ast})^{r-r_{x}}$, and $\delta^{i}(\psi_{\g}\C)_{x}=\overline{ \Delta^{i}(\U_{x},\nu_{x})}$.  Since the involution does not change $\W(\U_{x},\nu_{x})$, we then conclude that $\Su_{x}(\psi_{\g}\C)=\W^{unif}(\U_{x},\nu_{x}).$
 
As we will show later on in Proposition \ref{p6.1}, $\W(\U_{x},\nu_{x})$ is proper in $(\C^{\ast})^{r_{x}}$. Thus $\Su_{x}(\psi_{\g}\C)$ is proper in $(\C^{\ast})^{r}$. Finally, (\ref{3.3}) follows from Proposition \ref{p2.3} together with the fact that the involution keeps $\Delta^{i}(\U_{x},\nu_{x})$ unchanged.
\end{proof}

\br Proposition \ref{p3.1} is similar to, but more detailed than \cite[Lemma 3.18]{Bu}, which in turn is an extension of \cite[2.2.5]{Sab}. \cite{Bu} draws two conclusions out of this calculation: that one needs and can introduce $\W^{unif}$, as here, and $\Su^{unif}$. The last conclusion is not correct, as Proposition \ref{p3.1} points out, namely $\Su^{unif}$ should stay simply $\Su$ in \cite{Bu}. N. Budur \cite{Bu2} has communicated to us that, indeed, our Proposition \ref{p3.1} above is the correct version, and therefore it simplifies all the statements in \cite{Bu} concerning $\Su^{unif}$, by replacing it with $\Su$.
\er

Since $\Su_{x}(\psi_{\g}\C)$ is proper in $(\C^{\ast})^{r}$, we say that the Sabbah specialization complex is a torsion $\Gamma_{r}$-module complex, i.e., the stalks of its cohomology sheaves are torsion $\Gamma_{r}$-modules.

\bex \label{ex3.4} Nodal curve case:  Choose $g=x_{1}^{2}-x_{2}^{2}(x_{2}+1)$, and $x=(0,0)$. Then $H_{1}(\U_{x},\Z)=\Z^{2}$ and $\V(\U_{x})=\lbrace (1,1)\rbrace$. The map $\nu_{x}^{\ast}$ gives an embedding: $\C^{\ast}\to (\C^{\ast})^{2}$, where $\nu_{x}^{\ast}(t)=(t,t)$. So $\W(\U_{x},\nu_{x})=\lbrace 1 \rbrace $.
\eex

\bex \label{ex3.5} Choose $g_{1}=x_{1}^{2}-x_{2}^{2}(x_{2}+1)$, and let $g_{2}$ be a generic linear polynomial. Set $x=(0,0)$. Then $\U_{x}$ is homotopy equivalent to $S^{1}\times (S^{1}\bigvee S^{1})$, and $H_{1}(\U_{x},\Z)=\Z^{3}$. Hence $\V(\U_{x})=V(s_{1}s_{2}s_{3}=1)$. The map $\nu_{x}^{\ast}$ gives the embedding $(\C^{\ast})^{2}\to (\C^{\ast})^{3}$, where $\nu_{x}^{\ast}(t_{1},t_{2})=(t_{1},t_{1},t_{2})$. So $\W(\U_{x},\nu_{x})=V(t_{1}^{2}t_{2}=1) \subset(\C^{\ast})^{2} $.
\eex

Next, we give some vanishing results for $\Delta_{i}(\U_{x}, \nu_{x})$.
\bex \label{ex3.1} Assume that $D$ is a normal crossing divisor (NC) at $x$. Then $\U_{x}^{ab}$ is contractible.
Thus  $\V(\U_{x})=\lbrace \fst\rbrace$, and $\W(\U_{x},\nu_{x})=im(\nu_{x}^{\ast}) \cap \V(\U_{x})=\lbrace \fst \rbrace$. So $\Su_{x}(\psi_{\g}\C)= \W^{unif}(\U_{x},\nu_{x})=(\C^{\ast})^{r-r_{x}}$. If $r_{x}\geq 2$, then $\Delta_{i}(\U_{x},\nu_{x})=1$ for all $i$.
\eex

\bex \label{ex3.2} Assume that $D$ is an isolated non-normal crossing divisor at $x$. Then $\U_{x}^{ab}$ has the homotopy type of a bouquet of ($n-1$)-spheres, see \cite[Theorem 3.2]{DL}. So $A_{0}(\U_{x})=\Gamma_{r_{x}}/m_{\fst}$, and $A_{i}(\U_{x})=0$ for $0<i<n-1$. It follows that $\Delta_{i}(\U_{x},\nu_{x})=1$ for all $i<n-1$,  if $r_{x}\geq 2$.
\eex

\bex \label{ex3.3} Choose $x\in \bigcap_{i=1}^{r} D_{i}$, where $r\geq 2$. Assume that $D_{1}$ is smooth and $D_{1}$ intersects the other components of $D$ transversally in $B_{x}$. 
Set $\U_{x}^{\prime}:= B_{x}\setminus \bigcup_{i=2}^{r} D_{i}$. If follows that $\U_{x}^{ab} \simeq (\U_{x}^{\prime})^{ab}\times \mathbb{R} $ (see \cite[Remark 3.3]{DM}), hence $A_{i}(\U_{x})=A_{i}(\U_{x}^{\prime})\otimes_{\C} \dfrac{\C[t_{1},t_{1}^{-1}]}{\langle t_{1}-1\rangle}  $. Then  $\W(\U_{x},\nu_{x})$ is proper in the sub-torus of $(\C^{\ast})^{r}$ defined by $\lbrace t_{1}=1\rbrace$, so it has at least codimension $2$ in $(\C^{\ast})^{r}$, and $\Delta_{i}(\U_{x},\nu_{x})=1$ for all $i$.
\eex


\section{Divisibility results for hypersurfaces transversal at infinity}\label{S4}
From now on, we always assume that the projective hypersurface $V$ is transversal (in the stratified sense) to the hyperplane at infinity $H$ (i.e., $g=\prod_{i=1}^{r} g_{i}$ is {\it transversal (or regular) at infinity}). Then the affine hypersurface $D=\{g=0\}$ is homotopy equivalent to a bouquet of $(n-1)$-spheres, i.e.,  
\begin{equation}\label{4.1} D\sim \bigvee_{\mu} S^{n-1},\end{equation} where $\mu$ denotes the number of spheres in the above join (cf. \cite[page 476]{DP}). It is shown in loc.cit. that $\mu$ can be determined topologically as the degree of the gradient map associated to $f$.

In this section, we give  a general divisibility result (compare also with \cite[Theorem 3.2]{DM}) for the Alexander varieties of $\U=\mathbb{CP}^n\setminus (V \cup H)$ by using the associated {\it residue complex}.  


\subsection{Residue complex and its stalks}\label{s4.1} 

First, let us introduce the residue complex.
\bd Let $j:\U \hookrightarrow \CP$ be the inclusion map, and let $\K$ be the $\Gamma_r$-local system on $\U$ defined in Section \ref{s3.1}. The {\it residue complex $\R$ associated to $\U$} is defined by the distinguished triangle: \be \label{4.2}
j_{!}\K \to Rj_{\ast}\K \to \R \overset{[1]}{\to}
\ee
\ed

\br The residue complex is a generalization of the peripheral complex used in \cite{Max}.
In fact, when $r=1$, the residue complex $\R$ as defined here equals the (shifted) peripheral complex $\R[-2n]$  as defined by Cappell and Shaneson, see \cite{CS} or \cite{Max}.
\er

Recall from Section \ref{s3.1} that $\mathcal{D} \K  = \overline{\K}[2n]$. Then we have: \begin{center}
$\D j_{!}\K=  Rj_{\ast}\overline{\K}[2n]$, and $\D Rj_{\ast}\K= j_{!}\overline{\K}[2n]$.
\end{center} 
So, by dualizing (\ref{4.2}), we see that, up to a shift,  $\R$ is a self-dual (i.e., $\R \cong \mathcal{D} \overline{\R}[-2n+1] $), $\Gamma_r$-perverse sheaf. Moreover, by comparing the two distinguished triangles (\ref{3.1}) and (\ref{4.2}), we have
that $$\R\vert_{D}=\psi_{\g}\C.$$ Then similar properties hold for $\psi_{\g}\C$, since $D$ is open in $V\cup H$.

\medskip

Next, we compute the stalk cohomology of $\R$.  It is clear from definition that $\R$ has compact support on $V\cup H$. Fix a Whitney b-regular stratification $\ST$ of $V$. By the transversality assumption, there is a stratification $\ST^{\prime}$ of $V\cup H$ with strata of the form: $H\setminus V\cap H$, $S\cap H$ and $S\setminus S\cap H$, where $S\in \ST$.  Then $\Su_{x}\R$ is constant for $x$ in a given stratum of $\ST^{\prime}$.
The computation of stalk cohomology for $\R$ will be divided into the following 3 cases, according to the type of strata in $\ST^{\prime}$:

(1) If $x\in D=V\setminus H$, then $\R\vert_{D}=\psi_{\g}\C $, hence $\mathcal{H}^{i}(\R)_{x}$ for $x\in D$ is already computed in Section \ref{s3.3}.  In particular, we have in this case that $$\Su_{x}(\R)=\W^{unif}(\U_{x},\nu_{x}).$$

(2) If $x\in H \setminus V$, 
the link pair of $H \setminus V$ is  $(S^{1}, \emptyset)$, and the corresponding generator $\gamma_{\infty}$ is mapped to $\prod_{i=1}^{r} t_{i}^{-d_{i}}$ under the representation defining the local system $\mathcal{L}$.  So, in this case, 
\begin{center}
$\mathcal{H}^{i}(\R)_{x}\cong \left\{ \begin{array}{ll}
\Gamma_{r}/\langle \prod_{i=1}^{r} t_{i}^{-d_{i}}=1 \rangle, & i=1, \\
0, & \text{otherwise,}\\
\end{array}\right.$
\end{center} hence $$\Su_{x}(\R)=V(\prod_{i=1}^{r} t_{i}^{d_{i}}=1).$$

(3) If $x\in H \cap V$, assume that $x\in S\cap H$ with $\dim_{\C}S>0$, where $S\in \ST$.
Set $\U_{x}^{\prime}:= B_{x}-V$. The transversality assumption implies that $\U_{x}^{ab} \simeq (\U_{x}^{\prime})^{ab}\times \mathbb{R} $, and the deck group action on the product factor $\mathbb{R}$ is given as in the previous case. Choose an arbitrary point $x_{S}\in S\setminus S\cap H \subset D$. Then the K\"{u}nneth formula yields that 
\[
\begin{aligned}
 \mathcal{H}^{i}(\R)_{x} \cong 
 &   H^{i}(\U_{x_{S}}, \K\vert_{\U_{x_{S}}}) \otimes_{\Gamma_{r}}  \Gamma_{r}/\langle \prod_{i=1}^{r} t_{i}^{-d_{i}}=1 \rangle  \\  
& \oplus  H^{i-1}(\U_{x_{S}}, \K\vert_{\U_{x_{S}}}) \otimes_{\Gamma_{r}}  \Gamma_{r}/\langle \prod_{i=1}^{r} t_{i}^{-d_{i}}=1 \rangle .    
\end{aligned}
\]
In particular, Lemma \ref{l2.2} implies that in this case we get: 
$$\Su_{x}(\R) =  \W^{unif}(\U_{x_{S}}, \nu_{x_{S}}) \bigcap  V(\prod_{i=1}^{r} t_{i}^{d_{i}}=1).$$ 

The above stalk calculations, combined with Proposition \ref{p6.1} below, show that $\Su_{x}(\R)$ is always proper, for any $x \in V \cup H$. So $\R$ is a torsion $\Gamma_{r}$-module complex. Moreover, the involution does not change $\Su_{x}(\R)$ for any $x$.


\subsection{Divisibility results for Alexander varieties and characteristic varieties}
The following theorem is proved by using methods similar to those of \cite{Max}, where the second named author obtained general divisibility results for the classical (one-variable) Alexander invariants of the hypersurface complement $\U$ by making use the peripheral complex.

\bt \label{t4.1} Assume that the hypersurface $V \subset \CP$ is transversal to the hyperplane at infinity $H$. Fix a Whitney stratification $\ST$ of $V$, and let 
$\ST^{\prime}$ be the induced stratification of $V \cup H$ as above. Then, for any $0\leq k\leq n-1$, we have that 
\be  \label{4.3}
\bigcup_{i=0}^{k} \V_{i}(\U) =\bigcup_{i=0}^{k} \W_{i}(\U) \subset \bigcup_{S\subset D_{1}} \W^{unif}(\U_{x_{\ST}},\nu_{x_{S}}),
\ee
where the last union is over the strata $S$ of $\ST^{\prime}$ contained in a fixed irreducible component $D_{1}$ of $D$ with $\dim_{\C} S\geq n-k-1$ and $x_{S}\in S$ is arbitrary.   Moreover, 
\be  \label{4.4}
  \bigcup_{i=0}^{n-1} \V_{i}(\U) =\bigcup_{i=0}^{n-1} \W_{i}(\U)  \subset V(\prod_{i=1}^{r} t_{i}^{d_{i}}=1).
\ee
In particular,  $\W_{i}(\U)$ and $\V_{i}(\U)$ are proper in $(\C^{\ast})^{r}$ for all $i\leq n-1$. 
\et
\begin{proof} The proof is divided into 3 steps. 
\medskip 
\item[Step 1:] We first compute the Alexander modules $A_{i}(\U)$ ($i\leq n-1$) in terms of the residue complex $\R$.

Let $u$ and $v$ be the inclusions of $\CP\setminus V_{1}$ and respectively $V_{1}$ into $\CP$.  By applying the compactly supported hypercohomology  functor to the distinguished triangle $$v_{\ast}v^{!}j_{!}\K\to j_{!}\K \to u_{\ast}u^{\ast}j_{!}\K \overset{[1]}{\to},$$ we have the following long exact sequence:
 \be  \label{les}
\cdots \to H^{2n-i}(V_{1}, v^{!}j_{!}\K ) \to H^{2n-i}(\CP, j_{!}\K) \to H^{2n-i}(\CP\setminus V_{1}, u^{\ast} j_{!}\K)\to \cdots
\ee
Since $\U$ is smooth of complex dimension $n$, we have that $\K[n] \in \mbox{Perv}(\U) $. The inclusion $k:\U \hookrightarrow\CP\setminus V_{1}$ is a quasi-finite affine morphism, hence  $k_{!}\K[n] \in \mbox{Perv}(\CP\setminus V_{1})$ (e.g., see \cite{D2}, Corollary 5.2.17).
Since $\CP\setminus V_{1}$ is an $n$-dimensional affine variety, Artin's  vanishing theorem (e.g., see \cite{Sc}, Corollary 6.0.4) yields that:  
\be  \label{4.5}
H^{2n-i}(\CP\setminus V_{1}, k_{!}\K)=0 \text{ for }  i<n.
\ee
Note that $u^{\ast} j_{!}\K=k_{!}\K$, and recall that by Remark \ref{r3.1} we have that 
$$A_{i}(\U)\cong H^{2n-i}_{c}(\U, \K)\cong H^{2n-i}(\CP, j_{!}\K).$$ 
Therefore, the long exact sequence (\ref{les}), together with the vanishing from (\ref{4.5}), yields that: 
$$A_{i}(\U)\cong H^{2n-i}(V_{1}, v^{!}j_{!}\K) \text{ for } i<n-1,$$
and $A_{n-1}(\U)$ is a quotient of the $\Gamma_{r}$-module $H^{n+1}(V_{1}, v^{!}j_{!}\K)$.

\medskip

\item[Step 2:] Next, we prove the divisibility result (\ref{4.3}) for $k=n-1$.

 Let $c$ be the  map from $V_{1}$ to a point.  Then
\[
\begin{aligned}
H^{2n-i}(V_{1}, v^{!}j_{!}\K)
& \cong H^{2n-i}(Rc_{\ast} v^{!}j_{!}\K) \\
& \overset{(1)}{\cong} H^{-i}(Rc_{\ast}v^{!}j_{!}\D (\overline{\K})) \\
& {\cong} H^{-i}( Rc_{\ast} \D (v^{\ast}Rj_{\ast} \overline{\K})) \\
& \cong H^{-i}( \D(Rc_{!}(\overline{\R} \vert_{V_{1}})) )
\end{aligned}
\]
where (1) follows from $\mathcal{D} \overline{\K}  \cong \K[2n]$.

Note that $\D(Rc_{!}(\overline{\R} \vert_{V_{1}}))= R\Hom_{\Gamma_{r}}(Rc_{!}(\overline{\R} \vert_{V_{1}}),\Gamma_{r})$. Then we have the following spectral sequence  (e.g., see \cite[p.243]{B}) \be \label{ss4.1}
\Ext^{q}_{\Gamma_{r}}(H^{s}_{c}(V_{1},\overline{\R}) , \Gamma_{r}) \Longrightarrow H^{2n-q-s}(V_{1}, v^{!}j_{!}\K). 
\ee
On the other hand, $\D (\D \K)  =\K$, so we have a similar spectral sequence:\be \label{ss4.2}
\Ext^{q}_{\Gamma_{r}}(H^{2n-s}(V_{1}, v^{!}j_{!}\K), \Gamma_{r}) \Longrightarrow H^{s+q}_{c}(V_{1},\overline{\R} )
\ee
 By using the same arguments as in the proofs of Proposition \ref{p2.1}  and Proposition \ref{p2.2}, we obtain  
$$\bigcup_{i=0}^{2n} \Su(H^{2n-i}(V_{1}, v^{!}j_{!}\K))= \bigcup_{i=0}^{2n} \Su(H^{i}_{c}(V_{1},\overline{\R})).$$

Recall from Step $1$ above that $A_{i}(\U)\cong H^{2n-i}(V_{1}, v^{!}j_{!}\K)$ for $i< n-1$, and $A_{n-1}(\U)$ is a quotient of the $\Gamma_{r}$-module $H^{n+1}(V_{1}, v^{!}j_{!}\K)$. Then by Lemma \ref{l2.1} we have that $\W_{i}(\U) \subset \Su(H^{2n-i}(V_{1}, v^{!}j_{!}\K)$ for $i\leq n-1$. Thus:
\be 
\bigcup_{i=0}^{n-1} \W_{i}(\U) \subset  \bigcup_{i=0}^{n-1} \Su(H^{2n-i}(V_{1}, v^{!}j_{!}\K)) \subset \bigcup_{i=0}^{2n} \Su(H^{i}_{c}(V_{1},\overline{\R})). 
\ee
By using the compactly supported hypercohomology long exact sequence for the inclusion of strata of $V_{1}$, we get by Lemma \ref{l2.1}  that $$\Su (H_{c}^{i}(V_{1}, \overline{\R}))\subset \bigcup_{S\subset V_{1}} \Su (H_{c}^{i}(S,\overline{\R}\vert_{S})),$$ 
where $S$ runs over all strata of $\ST^{\prime}$ contained in $V_{1}$. Note that $H_{c}^{i}(S, \overline{\R}\vert_{S})$ is the abutment of a spectral sequence with the $E_{2}$-term defined by $E_{2}^{s,q}=H^{s}_{c}(S, \mathcal{H}^{q}(\overline{\R}))$, and $\Su^{i}_{x}(\R)$ is constant along the stratum $S$.  By using the exactness of localization and the corresponding localized spectral sequence, it follows that if $\lambda \notin \Su(\overline{\R}\vert_{S})$ then $\lambda \notin \Su (H_{c}^{i}(S,\overline{\R}\vert_{S}))$. So $$\Su (H_{c}^{i}(S,\overline{\R}\vert_{S})) \subset \Su(\overline{\R}\vert_{S})=\Su(\R\vert_{S}),$$
where the last identification follows from the fact that the involution does not change $\Su_{x}(\R)$, see the analysis of the stalks of $\R$ in Section \ref{s4.1}.
 
If $S\subset D_{1}$, then  it follows from the discussion of the previous section that $\Su(\R\vert_{S}) =\W^{unif}(\U_{x_{S}},\nu_{x_{S}})$, for $x_S \in S$ arbitrary. Similarly, if $S\subset V_{1}\cap H$,  there exists a stratum $S^{\prime}\subset D_{1}$ such that $\Su(\R\vert_{S}) =\W^{unif}(\U_{x_{S^{\prime}}},\nu_{x_{S^{\prime}}}) \cap V(\prod_{i=1}^{r} t_{i}^{d_{i}}=1) $, where $x_{S^{\prime}} \in S^{\prime}.$ Altogether,  \be 
\bigcup_{i=0}^{n-1} \W_{i}(\U) \subset   \bigcup_{S\subset D_{1}} \W^{unif}(\U_{x_{\ST}}, \nu_{x_{S}}).
 \ee
 
If we replace $V_{1}$ by the hyperplane at infinity $H$, then the stalk calculations for $\R$ over $H$ yield that \be \bigcup_{i=0}^{n-1} \W_{i}(\U) \subset V(\prod_{i=1}^{r} t_{i}^{d_{i}}=1).
\ee  
 \medskip
\item[Step 3:] We use the Lefschetz hyperplane section theorem to complete the proof of the theorem.

Fix $1\leq k <n-1$. Consider $L=\mathbb{CP}^{k+1}$  a generic  $(k+1)$-dimensional  linear subspace of $\CP$, such that $L$ is transversal to $V\cup H$. Then $W=L\cap V$ is a $k$-dimensional, reduced hypersurface in $L$, which is transversal to the hyperplane at infinity $H\cap L$ of $L$. Moreover, by the transversality assumption, $L\cap D$ has a Whitney stratification induced from $D$, with strata of the form $S\cap L$.

By applying the Lefschetz hyperplane section theorem (\cite{D1}, (1.6.5)) to the section of $\U$ by $L$, we see that the inclusion $\U\cap L \hookrightarrow\U$ is a homotopy $(k+1)$-equivalence. Hence $\W_{i}(\U\cap L)=\W_{i}(\U)$ for $i\leq k$. Using the fact that the link pair of a stratum $S\cap L$ in $L\cap D$ is the same as the link pair of $S$ in $D$, the claim follows by reindexing (replace $s$ by $s-(n-k-1)$, where $s=\dim_{\C}S$).
\end{proof}

\br  \label{r4.2}
Similar results are stated in \cite{DM}  for $\W^{i}(\U)$ ($i\leq n-1$) , see Theorem 3.2 and Theorem 3.6 in loc.cit. However, by using (\ref{2.2}), it can be seen that the statement of Theorem \ref{t4.1} is sharper than the above mentioned results from \cite{DM}. For the curve case, (\ref{4.4}) is first proved by Libgober in \cite[Corollary 3.3]{L92}.
\er
\bc  $A_{n}(\U)=0$ if and only if $\mu=0$.

\ec
\begin{proof}
If $\mu=0$, then $\U$ admits an infinite cyclic cover which has the homotopy type of a finite CW-complex of dimension $(n-1)$ (e.g., the generic fibre of $g$, see \cite[Proposition 4.6(c)]{LiM}). It then follows that the universal abelian cover $\U^{ab}$ has the homotopy type of an $(n-1)$-dimensional CW-complex, so $A_{n}(\U)=0$.

Conversely, if $A_{n}(\U)=0$, then $\V_{n}(\U)=0$, so $\bigcup_{i=0}^{n-1}\V_{i}(\U)=\bigcup_{i=0}^{n}\V_{i}(\U)= \V(\U) $. Theorem \ref{t4.1} shows that $\V(\U)$ is proper in $(\C^{\ast})^{r}$, hence by Remark \ref{r2.3} we get that $\chi(\U)=0$. Thus $\mu=(-1)^{n}\chi(\U)=0$.
\end{proof}

\br \label{r4.3}  A typical example with $\mu=0$ is the affine hypersurface complement $M$ defined in the introduction. Moreover, in the isolated singularities case, the condition $\mu=0$ implies that $\U$ is the affine complement of a homogeneous polynomial, up to a change of coordinates (see \cite[Theorem 1]{Huh} or the reformulation from \cite[Theorem 4.3]{Liu}). The hypersurface complement $\U$ with the transversality assumption and $\mu=0$ exhibits many similar properties with $M$, e.g., $\V(\U)\subset V(\prod_{i=1}^{r} t_{i}^{d_{i}}=1) $. For more results in this direction, see \cite[Proposition 4.2]{Liu} and \cite[Proposition 4.6]{LiM}.
\er 

\bex  Let $\U$ be the complement of the arrangement of 4 lines in $\C^{2}$ defined by $g = x(x-y)(x+y)(2x-y+1)$. Then $g$ is transversal at infinity and the line $\lbrace 2x-y+1=0\rbrace$ intersects the other components transversally. So $\V_{1}(\U)\subset V(t_{4}=1) \cap V(t_{1}t_{2}t_{3}t_{4}=1)$. (In fact, these two sets are equal, see \cite[Examples 2.5 and 3.7]{Su00}.)
\eex


\section{A more general setting}\label{S5}

Let $V\subset \CP$ be, as before, a degree $d$ hypersurface defined by the homogeneous polynomial $f=\prod_{i=1}^{r} f_{i}$, and set $$M^{\ast}:=\CP \setminus V.$$ Then $f$ also defines a hypersurface $f^{-1}(0)$ in $\C^{n+1}$, the affine cone on $V$. The complement $$M:=\C^{n+1} \setminus f^{-1}(0)$$ is the total space of the global Milnor fibration $f:M \to \C^{\ast}$, whose fiber $$F:=f^{-1}(1)$$ is called the Milnor fiber of $f$. Then there exists a $d$-fold covering map $p_{d}:  F\to M^{\ast}$.

The Hopf map induces a fibration: 
$$p: M \to M^{\ast}.$$
As we will show later on in the proof of Proposition \ref{p5.2} below, the transversality assumption implies that  $\pi_{1}(\U)\cong \pi_{1}(M)=G$, so $\lbrace \gamma_{i} \rbrace_{1\leq i \leq r}$ are also generators of $H_{1}(M,\Z)\cong \Z^{r}$.

In this section, we pass on the divisibility results from $\U$ to $M$ and $M^{\ast}$, respectively. As an application, we obtain vanishing results for the homology Alexander polynomials $\Delta_{i}(M)$, for $i\leq n-1$.


\subsection{Divisibility results}\label{s5.1}
To obtain the relation between the characteristic varieties of $M$ and $M^{\ast}$, we need to use the {\it total turn monodromy operator} defined in \cite{CDO}.  The following description can be found in \cite[Section 6.4]{D2}.

Fix a base point $a\in M$ and denote by $\sigma _{a}$ the loop $t\mapsto \exp(2\pi i t)a$ for $t\in[0,1]$. Choosing a generic line $L$ passing through $a$, transversal to $f^{-1}(0)$, and close to the line $\C a$ (which contains the loop $\sigma_a$), we see that the element $\sigma_{a}\in \pi_{1}(M,a)$ is given by the product (in a certain order) of the elementary loops $\sigma_{j}$ for $j=1,\cdots,d$, based at $a$ and associated to the intersection points in $L\cap f^{-1}(0)$. 

By \cite[Proposition 4.1.3]{D1} we have that $$H_{1}(M^{\ast},\Z)=\Z^{r-1}\oplus \left(\Z/\gcd(d_{1},\cdots,d_{r})\Z\right). $$
Moreover, $H_{1}(M^{\ast},\Z)\cong H_{1}(M,\Z)/\langle [\sigma_{1}]+\cdots+[\sigma_{d}]\rangle$, and note that $$[\sigma_{a}]=[\sigma_{1}]+\cdots+[\sigma_{d}]= d_{1}[\gamma_{1}]+\cdots+d_{r}[\gamma_{r}]$$ in $H_{1}(M,\Z)$.

For $\lambda=(\lambda_{1},\cdots,\lambda_{r})\in (\C^{\ast})^{r}=\Hom(G,\C^{\ast})$, consider the corresponding rank one local system $\K_{\lambda}$ on $M$ defined by the representation 
$\rho_{\lambda}\in {\Hom}( \pi_{1}(M,a), \C^{\ast})$, which sends each generator $\gamma_{i}$ of $H_{1}(M,\Z)$ to $\lambda_{i}$. We define the {\it total turn monodromy operator} of $\K_{\lambda}$ to be the invertible operator:
$$T(\K_{\lambda})=\rho_{\lambda}(\sigma_{a})\in Aut(\C)\cong \C^{\ast}.$$ Since $[\sigma_{a}]= d_{1}[\gamma_{1}]+\cdots+d_{r}[\gamma_{r}]$ in $H_{1}(M,\Z)$, we get that $$T(\K_{\lambda})=\prod_{i=1}^{r} \lambda_{i}^{d_{i}}\in \C^{\ast}.$$ As we will see below, the complex number $\prod_{i=1}^{r} \lambda_{i}^{d_{i}}$ plays a key role in describing the local system $R^{q}p_{\ast}\K_{\lambda}$. In fact, as shown in \cite[Section 6.4]{D2}, for any $x\in M^{\ast}$, we have: 
\begin{center}
$\mathcal{H} ^{q}(Rp_{\ast}\K_{\lambda})_{x}=\left\{ \begin{array}{ll}
E^{0}=\ker (\prod_{i=1}^{r} \lambda_{i}^{d_{i}}-1), & q=0, \\
E^{1}=\mbox{coker} (\prod_{i=1}^{r} \lambda_{i}^{d_{i}}-1), & q=1, \\
0, & \text{otherwise.}\\
\end{array}\right.$
\end{center} 

\medskip

If $\lambda\notin V(\prod_{i=1}^{r} t_{i}^{d_{i}}=1)$, i.e., $\prod_{i=1}^{r} \lambda_{i}^{d_{i}}\neq 1$, then $R^{q}p_{\ast}\K_{\lambda}=0$ for all $q$. Thus, for all $i$, we obtain in this case that $H^{i}(M,\K_{\lambda})=H^i(M^*,Rp_*\K_{\lambda})=0$. Then the duality result (\ref{2.8}) yields that $H_{i}(M,\K_{\overline{ \lambda}})=0$, or equivalently $\overline{\lambda} \notin \V(M)$. Thus $\overline{\V(M)} \subseteq V(\prod_{i=1}^{r} t_{i}^{d_{i}}=1)$, which also implies that 
$\V(M)\subset V(\prod_{i=1}^{r} t_{i}^{d_{i}}=1)$, since the involution $\lambda \mapsto \overline{\lambda}:=\lambda^{-1}$ keeps $V(\prod_{i=1}^{r} t_{i}^{d_{i}}=1)$ unchanged.  The latter inclusion was already mentioned in Remark \ref{r4.3}.

\medskip

If $\lambda\in V(\prod_{i=1}^{r} t_{i}^{d_{i}}=1)$, then $E^{0} \cong \C\cong E^{1}$. By the arguments used in 
\cite[Proposition 6.4.3]{D2}, there exists a rank one local system $\K^{\ast}_{\lambda}$ on $M^{\ast}$ such that $\K_{\lambda}=p^{\ast}\K^{\ast}_{\lambda}$. In fact, the representation of  $\K_{\lambda}$ on $M$ factors through the representation of $\K_{\lambda}^{\ast}$ on $M^{\ast}$: \be \label{5.10}
\xymatrix{
\pi_{1}(M) \ar[rd]^{\rho_{\lambda}} \ar[r]^{p_{\ast}}  & \pi_{1}(M^{\ast})  \ar[d]^{\rho_{\lambda}^{\ast}} \\
&   \C^{\ast}
}
\ee where $p_{\ast}$ is the epimorphism induced by  the Hopf map $p$. Here $\rho_{\lambda}$ and $\rho_{\lambda}^{\ast}$ are the  representation homomorphisms defining  $\K_{\lambda}$ and $\K_{\lambda}^{\ast}$, respectively. (Note that the relation $\sum_{i=1}^r d_i\gamma_i=0$ in $H_1(M^*,\Z)$ is sent by $\rho_{\lambda}^{\ast}$ to $\prod_{i=1}^{r} t_{i}^{d_{i}}=1$.)
The sheaves $R^{q}p_{\ast}\K_{\lambda}$  are local systems on $M^{\ast}$ for $q=0,1$, and trivial otherwise. As shown in \cite[Proposition 6.4.3]{D2}, the following construction gives the  representations defining $R^{q}p_{\ast}\K_{\lambda}$ ($q=0,1$): for $v \in E^{q}$, and $\alpha \in \pi_{1}(M^{\ast})$,  we set \begin{center}
$\alpha \cdot v = \rho_{\lambda}(\beta)(v)$
\end{center} where $\beta$ is any lifting of $\alpha$ under the epimorphism  $p_{\ast}$.  According to the commutative diagram (\ref{5.10}), this shows that the local systems $R^{q}p_{\ast}\K_{\lambda}$ ($q=0,1$) coincide  with $\K^{\ast}_{\lambda}$. 
As a consequence, the $E_2$-page of  the Leray  spectral sequence induced by the Hopf map $p:M \to M^*$, \be \label{ss5.1}
E_{2}^{s,q}=H^{s}(M^{\ast}, R^{q}p_{\ast}\K_{\lambda}) \Longrightarrow H^{s+q}(M,\K_{\lambda})
\ee reduces to $E_{2}^{s,q}=H^{s}(M^{\ast}, \K^{\ast}_{\lambda})$ for $q=0,1$ and trivial otherwise. Hence this spectral sequence degenerates at $E_{3}$.

\bp \label{p5.1} Under the above assumptions and notations, we have that $$\V_{i}(M)\subset  \V_{i}(M^{\ast}) \cup \V_{i-1}(M^{\ast}).$$ These two sets are equal if, in addition, one of the irreducible components of $V$ is a hyperplane.
\ep
\begin{proof}  By Remark \ref{r2.4}, it suffices to prove the claim for $\V^{i}(M)$.  Assume $\lambda \notin \V^{i}(M^{\ast}) \cup \V^{i-1}(M^{\ast})$. Then $E_{2}^{i,0}=H^{i}(M^{\ast}, \K^{\ast}_{\lambda})=0$ and $E_{2}^{i-1,1}=H^{i-1}(M^{\ast}, \K^{\ast}_{\lambda})=0$. The spectral sequence (\ref{ss5.1}) yields that $H^{i}(M, \K_{\lambda})=0 $, hence $\lambda \notin \V^{i}(M)$.

If, in addition, one of the irreducible components of $V$ is a hyperplane, then $M=M^{\ast}\times \C^{\ast}$ (see \cite[Proposition 6.4.1]{D2}). So the spectral sequence (\ref{ss5.1}) degenerates at $E_{2}$, and $H^{i}(M,\K_{\lambda})\cong H^{i}(M^{\ast}, \K^{\ast}_{\lambda}) \oplus H^{i-1}(M^{\ast}, \K^{\ast}_{\lambda})$. The claim follows.
\end{proof}

The next result establishes a close relation between the characteristic varieties of $M$ and $M^{\ast}$, respectively.
\bt \label{t5.1} For any $k\geq 0$, \be \label{5.1}
\bigcup_{i=0}^{k} \V_{i}(M)= \bigcup_{i=0}^{k} \V_{i}(M^{\ast}) .
\ee
\et
\begin{proof}
$(\subseteq)$ It follows from Proposition \ref{p5.1}.

$(\supseteq)$ The proof is done by induction on $k$.   By Remark \ref{r2.4}, it suffices  to prove the claim for $\V^{i}$.   For $k=0$, $\V^{0}(M)=\lbrace (\fst) \rbrace=\V^{0}(M^{\ast})$. Suppose $\lambda\in\bigcup_{i=0}^{k} \V^{i}(M^{\ast}) $. We may assume that $\lambda\notin\bigcup_{i=0}^{k-1} \V^{i}(M^{\ast}) $, and $\lambda \in \V^{k}(M^{\ast})$, for otherwise the proof is complete, by the induction hypothesis. These assumptions translate into $E_{2}^{s,q}=0$ for $s<k$
and $H^{k}(M^{\ast}, \K_{\lambda}^{\ast})\neq 0$.  Thus, \begin{center}
$H^{k}(M,\K_{\lambda})=E_{3}^{k,0}=E_{2}^{k,0}=H^{k}(M^{\ast}, \K_{\lambda}^{\ast}) \neq 0$.
\end{center}
So $\lambda \in \V^{k}(M)$ and the proof is completed.
\end{proof}

Set $G^{\ast}=\pi_{1}(M^{\ast})$ and $d_{i}^{\prime}=\dfrac{d_{i}}{\gcd(d_{1},\cdots,d_{r})}$.  The Hopf map induces an injective  map $ \Hom(\pi_{1}(M^{\ast}),\C) \to \Hom(\pi_{1}(M),\C)=(\C^{\ast})^{r}$. So we can regard $\Hom(\pi_{1}(M^{\ast}),\C)$ as a subset of $(\C^{\ast})^{r}$. In fact, we have that \begin{center}
${\Hom}(G^{\ast},\C^{\ast})=V(\prod_{i=1}^{r} t_{i}^{d_{i}}=1)$, and ${\Hom}(G^{\ast},\C^{\ast})^{0}=V(\prod_{i=1}^{r} t_{i}^{d_{i}^{\prime}}=1)$.
\end{center}  
The following proposition is a direct consequence of Theorem \ref{t5.1} and \cite[Corollary 6.5]{Liu}. In particular, it reproves (\ref{4.4}).
\bp  \label{p5.2} Assume that the hyperplane at infinity $H$ is transversal to the hypersurface $V\subset \CP$. Then, for any $k\leq n-1$, we have:
\be \label{5.2}
\bigcup_{i=0}^{k} \V_{i}(\U)=\bigcup_{i=0}^{k} \V_{i}(M)= \bigcup_{i=0}^{k} \V_{i}(M^{\ast}) \subset V(\prod_{i=1}^{r} t_{i}^{d_{i}}=1).
\ee
\ep
\begin{proof}
Since $H$ is transversal to $V$, it follows from  \cite[Corollary 6.5]{Liu}  that the affine hyperplane $\lbrace x_{0}=1 \rbrace$ in $\C^{n+1}$ is transversal to the affine hypersurface $f^{-1}(0)$. Note that $\U$ is isomorphic to $\lbrace x_{0}=1\rbrace\cap M$. The Lefschetz hyperplane section theorem 
shows that the natural inclusion $\U\cong \lbrace x_{0}=1\rbrace\cap M\hookrightarrow M$ is an $n$-homotopy equivalence; in particular, $\pi_{1}(\U)=\pi_{1}(M)=G$. Then,  for any $i\leq n-1$, we have: \begin{center}
$\W_{i}(\U)=\W_{i}(M)$, $\Delta_{i}(\U)=\Delta_{i}(M)$ and $\V_{i}(\U)=\V_{i}(M)$.
\end{center}
The claim follows now from Theorem \ref{t5.1}.
\end{proof}

Fix a  Whitney b-regular stratification $\ST$ of $V$. Since $V$ is transversal to $H$,  for any $S\in \ST$ we can find a point $x\in S\cap D$ such that $\U_{x}=M^{\ast}_{x}$, where $M^{\ast}_{x}=M^{\ast} \cap B_{x}$. Here $B_{x}$ is a small ball in $\CP$ centered at $x$. By using Proposition \ref{p5.2}, the divisibility results for $\U$ (cf. Theorem \ref{t4.1}) can be now translated into similar properties for  $M$ and $M^{\ast}$, respectively.  
\bt  \label{t5.2}  With the above assumption and notations, for any $k\leq n-1$ we have that:  \be  \label{5.3}
\bigcup_{i=0}^{k} \W_{i}(M)=\bigcup_{i=0}^{k} \V_{i}(M)  =  \bigcup_{i=0}^{k}  \V_{i}(M^{\ast})  \subset ( \bigcup_{S\subset V_{1}} \W^{unif}(M^{\ast}_{x_{S}}, \nu_{x_{S}}) ) \cap  V(\prod_{i=1}^{r} t_{i}^{d_{i}}=1),
\ee
\be  \label{5.4}
\bigcup_{i=0}^{k} \W_{i}(M^{\ast}) \subset  (\bigcup_{S\subset V_{1}} \W^{unif}(M^{\ast}_{x_{S}}, \nu_{x_{S}})) \cap V(\prod_{i=1}^{r} t_{i}^{d^{\prime}_{i}}=1) .
\ee where the unions are over the strata $S\in \ST$ contained in $V_{1}$ with $\dim_{\C} S\geq n-k-1$ and $x_{S}\in S$ is arbitrary.  
\et

\br \label{r5.1}
 Theorem \ref{t5.2} is a generalization of Theorem \ref{t4.1}. Indeed, let us replace $V$ by $V\cup H$. Set  $d_{0}=1$, and let $t_{0}$ be the corresponding coordinate for $H$ in the torus defined by the fundamental group.  As shown in Example \ref{ex3.3}, the transversality assumption implies that, for any $x\in H$,
 $\V^{unif}(\U_{x}) \subset V(t_{0}=1)$. Note that $$ V(\prod_{i=0}^{r} t_{i}^{d_{i}}=1) \bigcap V(t_{0}=1) =  V(\prod_{i=1}^{r} t_{i}^{d_{i}}=1)\times  V(t_{0}=1).$$
 Then  the divisibility result (\ref{5.3}) along $H$ gives (\ref{4.4}), while (\ref{4.3}) can be obtained from the divisibility results (\ref{5.3}) for $V_{1}$ by intersecting with $V(t_{0}=1)$.

Since $\W(\U_{x},\nu_{x})= \V(\U_{x})\cap im(\nu_{x}^{\ast})$, the divisibility results (\ref{5.3}) can be restated as saying that the global characteristic varieties (except the top one $\V_{n}(M^{\ast})$) are contained in the union of the local characteristic varieties along one irreducible component.

For the one-variable case, the above divisibility results imply the corresponding ones for the classical Alexander invariants of the infinite cyclic cover. More results in this direction are given in next section.
\er

Since the divisibility results of Theorem \ref{t5.2} hold for each irreducible component of $V$, we  get better estimates by taking the intersection over all the irreducible components.

\bc \label{c5.1}   Assume $r\geq 2$ and $V$ satisfies one of the following assumptions:
\begin{enumerate}
\item[(a)] $V$ is a normal crossing (NC) divisor at any point of the components $V_{1}, \cdots, V_{m}$ of $V$, where $m< r$.
\item[(b)]  Every irreducible component $\lbrace V_{1}, \cdots, V_{m} \rbrace$, with $m<r$,  is smooth and intersects transversally all other components of $V$.
\end{enumerate} Then $\dim\V_{i}(M^{\ast})\leq r-m-1$ for all $i\leq n-1$. In particular, if $m=r-1$, then,  $\bigcup_{i=0}^{n-1} \V_{i}(M)=\bigcup_{i=0}^{n-1} \V_{i}(M^{\ast})$  are just finite sets of torsion points.
\ec 
\begin{proof} It suffices to prove the case $m=1$, whereas the case $m\geq 2$ follows by taking intersections. 
If $V$ is only a NC divisor along $V_{1}$, then Examples \ref{ex3.1} shows that the local Alexander varieties along $V_{1}$ are contained in the subtori  $V(t_{1}=1)$, i.e., $\W^{unif}(\U_{x},\nu_{x})\subset V(t_{1}=1)$. On the other hand, as shown in Example \ref{ex3.3}, the transversality assumption of $(b)$ yields the same fact that  $\W^{unif}(\U_{x},\nu_{x})\subset V(t_{1}=1)$. Then in both cases, 
$$\bigcup_{i=0}^{n-1} \V_{i}(M^{\ast}) \subset V(\prod_{i=1}^{r} t_{i}^{d_{i}}=1) \bigcap  V(t_{1}=1),$$ 
with dimension at most $(r-2$).  
\end{proof}

\bex \label{ex5.1} Assume that $V$ is at most a NC divisor at any point $x\in V$. Then $$\bigcup_{i=0}^{n-1} \V_{i}(M)=\bigcup_{i=0}^{n-1} \V_{i}(M^{\ast}) =\lbrace \fst \rbrace.$$
\eex

\bex \label{ex5.2} Assume that $r=2$, and say that $V_{1}$ and $V_{2}$ have at most isolated singularities, and they intersect transversally. Then, for $x\in V_{1}\cap V_{2}$, $\W(M^{\ast}_{x})=\lbrace (1,1)\rbrace$; while for $x\in V_{i}$ ($i=1,2$), but $x\notin V_{1}\cap V_{2}$, $\W^{unif}(M^{\ast}_{x},\nu_{x})\subset \bigcup V(t_{i}=\lambda_{i})$, where $\lambda_{i}$ are the eigenvalues of the local Milnor fibration at the point $x\in V_{i}$.   So, by taking intersections, the divisibility result (\ref{5.3}) yields that
$$\bigcup_{i=0}^{n-1} \V_{i}(M^{\ast})\subset \lbrace (\lambda_{1}, \lambda_{2})\mid \lambda_{1}^{d_{1}}\lambda_{2}^{d_{2}}=1 \rbrace \cup \lbrace (1,1)\rbrace $$
 with $\lambda_{1}$, $\lambda_{2}$ being eigenvalues of the local Milnor fibration at the singular point of $V_{1}$, $V_2$, respectively.
\eex

Finally, it is  also natural to ask what can be said about $\V(M)=\V(M^{\ast})$. The next result follows at once from Remark \ref{r2.3} and Theorem \ref{t5.1}.

\bp \label{p5.3}   For $M$ and $M^{\ast}$ defined as above, we have that
 $$\W(M)=\V(M)=\V(M^{\ast})\subset V(\prod_{i=1}^{r} t_{i}^{d_{i}}=1).$$
If, in addition, $\chi(M^{\ast})\neq 0$, the above four sets are equal, and $\W(M^{\ast})=V(\prod_{i=1}^{r} t_{i}^{d^{\prime}_{i}}=1)$.
\ep


\subsection{Vanishing of Alexander polynomials}\label{s5.2} 
In this section, we derive vanishing results for the homology Alexander polynomials $\Delta_{i}(M)$, for $i\leq n-1$.

\bt \label{t5.3} Assume that the hypersurface $V\subset \CP$ is transversal to the hyperplane at infinity $H$. Then, for fixed $i\leq n-1$, the prime factors of $\Delta_{i}(\U)=\Delta_{i}(M)$ are among the prime factors of the local Alexander polynomials $\Delta_{q}(\U_{x_{S}},\nu_{x_{S}})$ associated to strata $S\subset \bigcap_{i=1}^{r}D_{i}$, with $\dim_{\C}S=  s\geq n-k-1$ and $0\leq q \leq n-s-1$.
\et
\begin{proof} The prime factors of $\Delta_{i}(\U)$ are in one-to-one correspondence to the codimension one irreducible hypersurfaces in $(\C^{\ast})^{r}$ contained in $\W_{i}(\U)$. Theorem  \ref{t4.1} shows that the prime factors of $\Delta_{i}(\U)$ $(i\leq n-1)$ are among the primes factors of the local Alexander polynomials $\Delta_{q}(\U_{x_{S}},\nu_{x_{S}})$ associated to strata $S\subset D$. We only need to show that  the local Alexander polynomials $\Delta_{q}(\U_{x_{S}},\nu_{x_{S}})$, with $S\nsubseteq \bigcap_{i=1}^{r}D_{i}$, are coprime with $\Delta_{i}(\U)$.
 
Assume that $x\in D$, but $x\notin \bigcap_{i=1}^{r}D_{i}$ (i.e., $r_{x}< r$). Then $\Delta_{q}(\U_{x},\nu_{x})$ is a polynomial in the variables $\lbrace t_{i} \rbrace_{i\in I_{x}}$. Since $r_{x}< r$, it follows that $\gcd (\Delta_{q}(\U_{x},\nu_{x}),\prod_{i=1}^{r} t_{i}^{d_{i}}-1)=1$.
On the other hand, by (\ref{4.4}) we get that $\Delta_{i}(\U)$ ($i\leq n-1$) divides $(\prod_{i=1}^{r} t_{i}^{d_{i}}-1)$. Altogether, we get that $\gcd (\Delta_{q}(\U_{x},\nu_{x}),\Delta_{i}(\U))=1$.
\end{proof}


Set $\f=(f_{1},\cdots,f_{r})$. Then $\f$ defines a polynomial map from $\C^{n+1}$ to $\C^{r}$. The following proposition is proved in \cite[Corollary 3.9]{Bu} by using Sabbah's generalization of A'Campo's formula (\cite[2.6.2]{Sab}).  
\bp \label{p5.4} $Z(\psi_{\f}\C)_{0}=(\prod_{i=1}^{r} t_{i}^{d_{i}}-1)^{-\chi(M^{\ast})}$, where $0$ is the origin of $\C^{n+1}$. 
\ep 

Since $f$ is homogeneous, the local Milnor fibration at the origin is fiber diffeomorphic to the global  Milnor fibration of $M$, see \cite[Exercise 3.1.13]{D1}. In particular, the corresponding local and global Alexander polynomials coincide.  
Then the next result follows at once from  Theorem \ref{t5.3}, Proposition \ref{p5.4} and Examples \ref{ex3.1}, \ref{ex3.3}.

\bc \label{c5.4} Assume that $V$ is transversal to $H$, $r\geq 2$, and $V$ satisfies one of the following assumptions:
\begin{enumerate}
\item[(1)]  $V$ is an essential hypersurface, i.e.,  $\bigcap_{i=1}^{r} V_{i}= \emptyset$.
\item[(2)]  One of the irreducible components of $V$, say $V_{1}$, is smooth and transversal to $\bigcap_{i=2}^{r} V_{i}$.
\item[(3)] $V$ is a normal crossing divisor along any point of $\bigcap_{i=1}^{r} V_{i}$.
\end{enumerate} Then $\Delta_{i}(\U)=\Delta_{i}(M)=1$ for all $i\leq n-1$, and 
$$\dim (\bigcup_{i=0}^{n-1} \V_{i}(\U))=\dim (\bigcup_{i=0}^{n-1}\V_{i}(M))=\dim (\bigcup_{i=0}^{n-1}\V_{i}(M^{\ast})) \leq r-2.$$ 
Furthermore, $$\Delta_{n}(M)=\big(\prod_{i=1}^{r} t_{i}^{d_{i}}-1\big)^{(-1)^{n}\chi(M^{\ast})}.$$  In particular,  
\be  \label{5.5}
(-1)^{n}\chi(F)=(-1)^{n} d \cdot\chi(M^{\ast})\geq 0.
\ee 
\ec
\begin{proof}
Since  $\W(M) \subset V(\prod_{i=1}^{r} t_{i}^{d_{i}}=1)$, then $\W(M)$ is proper in $(\C^{\ast})^{r}$. By using Proposition \ref{p3.1}, we have that $\delta^{i}(\psi_{\f}\C)_{0}=\Delta^{i}(M)=\Delta_{i-1}(M)$ for all $i$. Therefore,  
$$\big(\prod_{i=1}^{r} t_{i}^{d_{i}}-1\big)^{-\chi(M^{\ast})}=Z(\psi_{\f}\C)_{0}=\prod_{i} \delta^{i}(\psi_{\f}\C)_{0}^{(-1)^{i}}=\prod_{i} \Delta_{i}(M)^{(-1)^{i+1}}.$$
By using  Theorem \ref{t5.3} and Examples \ref{ex3.1}, \ref{ex3.3}, we see that in all three situations above we get: $\Delta_{i}(\U)=\Delta_{i}(M)=1$ for all $i\leq n-1$. Hence $\big(\prod_{i=1}^{r} t_{i}^{d_{i}}-1\big)^{-\chi(M^{\ast})}=\Delta_{n}(M)^{(-1)^{n+1}}$, so 
\be  \label{5.6}
  \Delta_{n}(M)=\big(\prod_{i=1}^{r} t_{i}^{d_{i}}-1\big)^{(-1)^{n}\chi(M^{\ast})}.
\ee
Since the degree of $\Delta_{n}(M)$ is non-negative, we must then have that 
\be  \label{5.7}
   (-1)^{n} \chi(M^{\ast})\geq 0.
\ee
\end{proof}

\br The inequality (\ref{5.7}) generalizes several known results obtained by different techniques:
\begin{itemize}
\item[(1)] The first situation generalizes \cite[Corollary 2.2]{D3} for essential hyperplane arrangements, where (\ref{5.7}) is proved by using $M_{0}$-tame polynomials. In fact, in this case $D$ is homotopy equivalent to a bouquet of $(n-1)$-spheres.   The same result is also proved in \cite[Proposition 2.1]{DJL} by induction.
\item[(2)] A typical example for the second situation is an affine hypersurface complement $\U$ with the transversality assumption at infinity. In fact, in this case we have that $(-1)^{n}\chi(\U)=\mu \geq 0$.
\end{itemize}
\er
A different proof of (\ref{5.7}) will be given in Corollary \ref{c6.1} of the next section.


\subsection{Hyperplane arrangements}\label{s5.3}
In this section, we focus on the special case of hyperplane arrangements.
Without loss of generality, we assume that $V$ is an essential hyperplane arrangement.
Note that $H_{1}(M^{\ast}_{x},\Z)=\Z^{r_{x}}$, so the first map in (\ref{3.2}) is an isomorphism. Moreover, since the local defining polynomial is homogeneous, we also have that: $$\V(M^{\ast}_{x})=\W(M^{\ast}_{x}) \subset V(\prod_{i\in I_{x}} t_{i}=1).$$  Then the following divisibility result for hyperplane arrangements holds:

\bt  \label{t5.4}  Let $M^{\ast}=\CP \setminus V$ be the complement of an essential hyperplane arrangement. Then, for any $k\leq n-1$,  \be   \label{5.9}
\bigcup_{i=0}^{k} \V_{i}(M)  =  \bigcup_{i=0}^{k}  \V_{i}(M^{\ast})  \subset  \bigcup_{S\subset V_{1}} \lbrace\V(M^{\ast}_{x_{S}}) \times V(\prod_{i \notin I_{x_S}} t_{i}=1)\rbrace,
\ee where the union is over the strata $S$ of $V_{1}$ with $\dim_{\C} S\geq n-k-1$, and $x_{S}\in S$ is arbitrary.  
\et
\begin{proof}
The claim follows from Theorem \ref{t5.2} and the observation that for any $x \in V$ we have: $$V(\prod_{i\in I_{x}} t_{i}=1) \cap V(\prod_{i=1}^{r} t_{i}=1) =V(\prod_{i\in I_{x}} t_{i}=1) \times V(\prod_{i \notin I_{x}} t_{i}=1) .$$
\end{proof}

\br Theorem  \ref{t5.4} can be used to derive nonresonance-type results for rank-one local systems, such as \cite[Theorem 6.4.18]{D2}.
\er

For $\lambda \in V(\prod_{i=1}^{r}t_{i}=1)$,  let $\K_{\lambda}^{\ast}$ be the corresponding local system on $M^{\ast}$.  
\bd $\lambda$ is called {\it generic at $x$}  if $\prod_{i \in I_{x}}\lambda_{i}\neq 1$ or $\prod_{i\notin I_{x}}\lambda_{i}\neq 1$. 
\ed
The next result should be compared with \cite[Theorem 5.3]{DJL}. 

\bc \label{c5.5} Let $M^{\ast}$ be the complement of an essential hyperplane arrangements. Assume that $\lambda$ is generic for any point $x\in V_{1}$. Then  \begin{center}
$H_{i}(M^{\ast},\K_{\lambda}^{\ast})=\left\{ \begin{array}{ll}
\C^{(-1)^{n} \chi(M^{\ast}) }, & i=n, \\
0, & \text{otherwise.}\\
\end{array}\right.$
\end{center} 
\ec
\begin{proof}
By Theorem  \ref{t5.4}, we see that if $\lambda$ is generic along $V_{1}$, then $\lambda \notin \bigcup_{i=0}^{n-1}\V_{i}(M^{\ast})$, and the claim follows.
\end{proof}


\section{Dwyer-Fried invariants and classical Alexander modules}\label{S6}
In this section, by using the Dwyer-Fried invariants, we recast some old and obtain new finiteness and divisibility results for the classical (infinite cyclic) Alexander modules of hypersurface complements.

\subsection{Dwyer-Fried invariants}\label{s6.1}
In this subsection, we introduce the {\it Dwyer-Fried sets}, and establish a relation with the characteristic varieties. As before, we assume that $X$ is a finite connected $n$-dimensional CW complex, with $\pi_{1}(X)=G$. Let $H_{1}(X,\Z)/{(\mbox{torsion})} \cong \Z^{r}$ be the maximal torsion free abelian quotient of $G$, where $r=b_{1}(X)$. 
The following brief introduction on the Dwyer-Fried sets can be found in \cite{Su12B}. 

Fix an integer $1\leq m \leq r$, and consider the regular covers of $X$, with deck groups $\Z^{m}$. Each such cover, $X^{\nu}\to X$, is determined by an epimorphism $\nu: H_{1}(X,\Z)\twoheadrightarrow \Z^{m} $. The induced  homomorphism in rational cohomology, $\nu^{\sharp}: \Q^{m}\hookrightarrow H^{1}(X,\Q)$, defines an $m$-dimensional $\Q$-vector subspace, $P_{\nu}=im(\nu^{\sharp})$, in $H^{1}(X,\Q)\cong \Q^{r}$. On the other hand, any $m$-dimensional $\Q$-vector subspace can be realized by some epimorphism  $ H_{1}(X,\Z)\twoheadrightarrow \Z^{m} $.

\bp  \cite{DF} The connected, regular covers of  $X$, with deck groups $\Z^{m}$, are parametrized by the Grassmannian of the $m$-planes in $H^{1}(X,\Q)$, via the correspondence:
\begin{align*}
\lbrace \text{regular } \Z^{m}\text{-covers   of } X  \rbrace  & \longleftrightarrow \lbrace \text{$m$-planes in } H^{1}(X,\Q) \rbrace \\
X^{\nu}\to X   & \longleftrightarrow  P_{\nu}=im(\nu^{\sharp}).
\end{align*}
\ep
  
\bd The {\it Dwyer-Fried invariants of $X$ are the subsets} 
$$\Omega^{k}_{m}(X)=\lbrace  P_{\nu}\in Gr_{m}(H^{1}(X,\Q)) \mid b_{i}(X^{\nu})<\infty \ \text{for} \ i\leq k \rbrace.$$
\ed

The $\Omega$-sets are homotopy-type invariants of $X$, see \cite[Lemma 3.5]{Su14}. The next result was first proved in \cite[Theorem 1]{DF}, and further improved in \cite{PS10}.

\bt \label{t6.1}  \cite[Corollary 6.2]{PS10} Assume that $X$ is a finite connected CW complex. Let $\exp$ denote the exponential map from $\C^{r}$ to $(\C^{\ast})^{r}$.  Then, for any $k\geq 0$, we have that
\be 
\Omega ^{k}_{m}(X)=\lbrace P_{\nu}\in Gr_{m}(H^{1}(X,\Q))  \mid \left(\exp(P_{\nu}\otimes \C)\cap  (\bigcup_{i=0}^{k} \V_{i}(X))\right) \text{ is finite} \rbrace.
\ee
\et

\bex  In the setup of Example \ref{ex5.2}, we get that $\bigcup_{i=0}^{n-1}\W_{i}(\U)$ is finite. So $\Omega_{2}^{n-1}(\U)\neq \emptyset$.
\eex


\subsection{Classical Alexander modules}\label{s6.2}
We now analyze in detail one  particular regular cover for the hypersurface complement $\U$ in $\CN$. In this subsection, we will not assume that the hypersurface $V\subset \CP$ is transversal to the hyperplane at infinity $H$.  

Recall that $H_{1}(\U,\Z)=\Z^{r}$. Let $lk$ denote the epimorphism $H_{1}(\U,\Z)\twoheadrightarrow \Z$, which maps the generators $\gamma_{i}$ to $1_{\Z}$. (Here $P_{lk}=\langle[1,\cdots,1]\rangle \in Gr_{1}( H^{1}(\U,\Q))$). Consider the  infinite cyclic cover $\U^{c}$ of $\U$ defined by $lk$, which is usually referred to as  
 the linking number infinite cyclic cover. The Alexander modules $H_{i}(\U^{c})$  are finitely generated $\Gamma=\C[t,t^{-1}]$-modules. If $H_{i}(\U^{c})$ is a torsion $\Gamma$-module, let $\omega_{i}(t)$ denote the corresponding Alexander polynomial of $H_{i}(\U^{c})$.  Many of the classical finiteness  and divisibility results about $H_{i}(\U^{c})$ for $i\leq n-1$ can be recovered by applying Theorems  \ref{t5.2} and \ref{t6.1} to this specific regular cover.
 
Note that if $V$ is irreducible (i.e., $r=1$), the infinite cyclic cover coincides with the universal abelian cover. In this case, $H_{i}(\U^{c})$ is a torsion $\Gamma$-module if and only if $\W_{i}(\U)$ is proper, and $\W_{i}(\U)$ consists of the collection of the roots of $\omega_{i}(t)$.

\medskip
 
For $x\in D=V \setminus H$, the composed map $H_{1}(\U_{x},\Z)\to H_{1}(\U,\Z)\overset{lk}{\to} \Z$ is an epimorphism, where the first map is induced by the natural inclusion $\U_{x}\to \U$. We will use the notation $lk_{x}$ for this composed map. Recall the notation  $g=\prod_{i=1}^{r}g_{i}$ introduced in Section 3.1. For $x\in D$, let $F_{x}$ be the local Milnor fibre of $g$ at $x$, which is homotopy equivalent to a finite CW complex. Then the corresponding local infinite cyclic cover $\U_{x}^{c}$ is homotopy equivalent to $F_{x}$, so $P_{lk_{x}} \in \Omega^{k} _{1}(\U_{x})$ for any $k\geq 0$.

\bp \label{p6.1}  For $x\in D$, $\V(\U_{x})$ is proper. Moreover, $\W(\U_{x},\nu_{x})$ is proper in $(\C^{\ast})^{r_{x}}$. 
\ep  
\begin{proof} The local infinite cyclic cover $\U_{x}^{c}$ is homotopy equivalent to the Milnor fibre $F_{x}$, so Theorem \ref{t6.1} yields  that $\exp(P_{lk_{x}} \otimes \C) \cap \V(\U_{x})$ is finite.  Note that $H_{1}(\U_{x},\Z)$ is torsion free.
If $\V(\U_{x})$ is not proper, then $\V(\U_{x})=(\C^{\ast})^{b_{1}(\U_{x})}$ and $\exp(P_{lk_{x}} \otimes \C) \cap \V(\U_{x})=\exp(P_{lk_{x}} \otimes \C)$ is infinite, which gives a contradiction.

For the second claim, note that the map $lk_{x}$ factors through $\Z^{r_{x}}$. In fact, we have a diagram 
\begin{center}
$\xymatrix{
\pi_{1}(\U_{x}) \ar[r] \ar@/^2pc/[rr]^{\nu_{x}} & H_{1}(\U_{x},\Z) \ar[rrd]^{lk_{x}} \ar[r] &  \Z^{r_{x}}   \ar[r] & H_{1}(\U,\Z)  \ar[d]^{lk} \\
& &  &  \Z 
}$
\end{center}
where the horizontal  map is the composed map in (\ref{3.2}). It follows that $\exp(P_{lk_{x}} \otimes \C) \subset im(\nu_{x}^{\ast})$. Recall that $\W(\U_{x},\nu_{x})= im(\nu_{x}^{\ast}) \cap \V(\U_{x})$. Then $$\exp(P_{lk_{x}} \otimes \C) \cap \W(\U_{x},\nu_{x})=\exp(P_{lk_{x}} \otimes \C) \cap \V(\U_{x})\cap im(\nu_{x}^{\ast})=\exp(P_{lk_{x}} \otimes \C) \cap \V(\U_{x}) .$$  So $\exp(P_{lk_{x}} \otimes \C) \cap \W(\U_{x},\nu_{x})$ is also finite. A similar argument as above then shows that $\W(\U_{x},\nu_{x})$ is proper in $(\C^{\ast})^{r_{x}}$.
\end{proof}

Set $\T=\lbrace (\alpha,\cdots,\alpha) \mid \alpha\in \C^{\ast}\rbrace \subset (\C^{\ast})^{r}$. It is clear that $\exp(P_{lk} \otimes\C)=\T$.
Similarly, set $\T_{x}=\lbrace (\alpha,\cdots,\alpha) \mid \alpha\in \C^{\ast}\rbrace \subset (\C^{\ast})^{r_{x}}$. 

\bp  \label{p6.2} Assume that $P_{lk} \in \Omega^{j}_{1}(\U)$. Then, for all $0\leq k\leq j$, $(\bigcup_{i=0}^{k} \W_{i}(\U)) \cap \T$ equals the collection of the roots of the Alexander polynomials $\lbrace\omega_{i}(t)\rbrace_{0\leq i \leq k}$. Similarly, for $x\in D$, $\W(\U_{x},\nu_{x}) \cap \T_{x}$ is the collection of all the eigenvalues for the local Milnor fibration of $g$ at $x$.
\ep 
\begin{proof}
If $(\alpha,\cdots, \alpha) \in \bigcup _{i=0}^{k}\W_{i}(\U)=\bigcup _{i=0}^{k}\V_{i}(\U)$, then  there exists $0\leq i\leq k$ such that $H_{i}(\U,\K_{\alpha})\neq 0$.
Since $P_{lk} \in \Omega^{j}_{1}(\U)$, the Milnor long exact sequence shows that  $\dim H_{i}(\U,\K_{\alpha})= N(i,\alpha)+N(i-1,\alpha)>0$, where $N(i,\alpha)$ is the number of the direct summands in the $(t-\alpha)$ part of $H_{i}(\U^{c})$, see \cite[Theorem 4.2]{DN}. So $\alpha$ must be a root of $\omega_{i}(t)$ or $\omega_{i-1}(t)$, and vice versa. 

Note that, for any $x\in D$, $P_{lk_{x}}  \in \Omega^{k} _{1}(\U_{x})$ for $k\geq 0$.  As shown in the proof of Proposition \ref{p6.1}, $\T_{x} \cap (\bigcup_{i=0}^{k} \W_{i}(\U_{x},\nu_{x}))=\exp(P_{lk_{x}} \otimes \C) \cap (\bigcup_{i=0}^{k} \V_{i}(\U_{x})) .$
So the second claim follows by a similar argument as the first part of the proof.   
\end{proof}

\br The statement of the above proposition should be known to experts, and it had been part of the mathematical folklore. For example, the corresponding claim for curve case can be found in \cite[Section 1.3.4]{L01}. We include it here since we were not able to find a written proof of it in the general case (compare also with the more recent result from \cite[Proposition 1.3]{BW15}).
\er

Note that, by combining  Theorems \ref{t5.2} and \ref{t6.1}, we get that, if $\bigcup \W^{unif}(\U_{x},\nu_{x}) \cap \T$ is finite along one of the irreducible components of $V\cup H$, then $P_{lk} \in \Omega_{1}^{n-1}(\U)$.

\medskip
 
Let us now fix some notations. Set $d_{0}=1$, and let $t_{0}$ be the corresponding coordinate for $H$ in the torus defined by the fundamental group. We analyze the local $\Omega$-sets in the following three cases:
\begin{enumerate}
\item[(a)] if $x\in D$, then $P_{lk_{x}} \in \Omega_{1}^{i}(\U_{x})$ for any $i\geq 0$, so $\W^{unif}(\U_{x},\nu_{x}) \cap \T =\W(\U_{x},\nu_{x}) \cap \T_{x}$ is always a finite set.  
\item[(b)] if $x\in H\setminus V\cap H$, then points of $V(\prod_{i=0}^{r} t_{i}^{d_{i}}=1) \cap \T$ satisfy $t_{0}=\alpha^{-d}$. Since $H\setminus V \cap H$ is smooth, then $\W^{unif}(\U_{x},\nu_{x}) =V(t_{0}=1)$, hence $\alpha^{d}=1$. Therefore, we have only finite choices for $\alpha$.
\item[({c})] The case when $x\in V \cap H$ will be treated in more detail below. More precisely, we aim to show that, under some good assumptions for $V \cap H$, we have that $P_{lk} \in \Omega_{1}^{n-1}(\U)$.
\end{enumerate}

\bp \label{p6.3}  Assume that $V$ and $H$ satisfy one of the following assumptions:
\begin{enumerate}
\item[(1)] $V \cup H$ is an essential hyperplane arrangement.
\item[(2)] $H$ intersects $V$ transversally.
\item[(3)]$V\cup H$ is a normal crossing divisor along any point of $V\cap H$.
 \end{enumerate}
Then $P_{lk} \in \Omega_{1}^{n-1}(\U)$, i.e., $H_{i}(\U^{c})$ is a torsion $\Gamma$-module for $i\leq n-1$.
\ep 

\begin{proof} In the first case,  $V\cup H$  is an essential hyperplane arrangement, so without loss of generality we can assume that, for $x\in V\cap H$, $x\in  V_{1}\cap \cdots \cap V_{k} \cap H$, where $1\leq k<r$. Then $\W^{unif}(\U_{x},\nu_{x}) \subset V(\prod _{i=0}^{k} t_{i}=1)$. Note that points of $V(\prod _{i=0}^{r} t_{i}=1)\cap V(\prod _{i=0}^{k} t_{i}=1)\cap \T$ satisfy $\alpha^{r-k}=1$ (as, in the hyperplane arrangements case, we have $d=r$). Here $\T\subset (\C^{\ast})^{r}$, where $(\C^{\ast})^{r}$ has coordinates $t=(t_{1},\cdots,t_{r})\in (\C^{\ast})^{r}$.  Since $k<r$, we can have only finitely many choices for $\alpha$.

In the second and third cases, the transversality and the normal crossing assumptions along $V\cap H$ both imply that $\W^{unif}(\U_{x},\nu_{x}) \subset V(t_{0}=1)$ for any $x\in V\cap H$, hence $\alpha^{d}=1$.
 \end{proof}

\br \label{r6.1}  Proposition \ref{p6.3} shows that in the above three cases the Alexander module $H_{i}(\U^{c})$ is a torsion $\Gamma$-module for all $i\leq n-1$. Note that by Theorem \ref{t5.2} we have that \begin{center}
$\bigcup_{i=0}^{k} \W_{i}(\U)  \cap \T   \subset \lbrace \bigcup_{S\subset V_{1}} \W^{unif}(\U_{x_{S}},\nu_{x_{S}}) \rbrace \cap  V(\prod_{i=0}^{r} t_{i}^{d_{i}}=1) \cap \T$,
\end{center} where the union is over the strata of one irreducible component $V_{1}$ of $V\cup V_{0}$ with $\dim_{\C} S\geq n-k-1$ and $x_{S}\in S$ is arbitrary. By using Proposition \ref{p6.2}, we obtain the corresponding divisibility result for $\omega_{i}(t)$ ($i\leq n-1$). \er

Let us now explain in more detail some consequences of our divisibility results: 
\begin{itemize}
\item[(1)] Assume that $V \cup H$ is an essential hyperplane arrangement. It was proved by Dimca that $H_{i}(\U^{c})$ is a torsion $\Gamma$-module for $i\leq n-1$, see \cite[Theorem 2.2]{D3}. To our knowledge, the following results are new: the divisibility result along $H$ shows that the roots of  $\omega_{i}(t)$ ($i\leq n-1$) have order $r!$; while the divisibility result along one of the irreducible components, say $V_{1}$, shows that the roots of $\omega_{i}(t)$  ($i\leq n-1$) are among  
$\lbrace \alpha \mid \alpha^{(r-1)!}=1 \rbrace$
and the zeros of local characteristic polynomial $\omega_{q}(x_{S})(t)$ associated to $x_{S}\in S\subset D_{1}$,  in the range $\dim_{\C} S=s\geq n-1-i$ and $0\leq q \leq n-1-s$.
 
\item[(2)]  Assume that $H$ intersects $V$ transversally. The transversality assumption implies that $\W^{unif}(\U_{x},\nu_{x}) \subset V(t_{0}=1)$ for any $x\in H$. The divisibility result along $H$ yields that the roots of $\omega_{i}(t)$ ($i\leq n-1$) have order $d$, a result already obtained by the second named author in \cite[Theorem 4.1]{Max}. Moreover, the divisibility result along $V_{1}$ shows that the roots of $\omega_{i}(t)$  ($i\leq n-1$) are among the roots of $\omega_{q}(x_{S})(t)$,  in the range $\dim_{\C} S=s\geq n-1-i$ and $0\leq q \leq n-1-s$. Here $\omega_{q}(x_{S})(t)$ is the $q$-th characteristic polynomial of the monodromy of the local Milnor fibration of $g$ at $x_{S}\in S\subset D_{1}$. So we also recover the second-named author's divisibility results for Alexander polynomials, see \cite[Theorem 4.2]{Max}.   
\end{itemize}

\br Assume that $V$ has only isolated singularities including at infinity. This is the case considered by Libgober in \cite{L94}. He showed that $P_{lk_{x}} \in \Omega_{1}^{i}(\U_{x})$ for any $i\geq 0$, if $x\in Sing_{\infty}(V)$. Then $P_{lk} \in \Omega_{1}^{n-1}(\U)$, i.e., $H_{i}(\U^{c})$ is a torsion $\Gamma$-module for $i\leq n-1$. It follows by the Lefschetz hyperplane section theorem that  $\omega_{i}(t)=1$ for $1\leq i\leq n-2$, see \cite{L94}, so the only interesting Alexander polynomial is $\omega_{n-1}(t)$. The divisibility result along $V$ yields that the roots of $\omega_{n-1}(t)$ (besides $1$) are among the roots of $\omega(x)(t)$, for $x\in Sing(V) \cup Sing_{\infty}(V) $.  Here, for $x\in D$, $\omega(x)(t)$ is the corresponding top characteristic polynomial of the local Milnor fibration of $g$ at $x$; while for $x\in Sing_{\infty}(V)$, $\omega(x)(t)$ is defined as in \cite[Definition 4.1]{L94}. We therefore recast Libgober's  divisibility result for the Alexander polynomials, see \cite[Theorem 4.3]{L94}. 
 \er

\bc \label{c6.1} Under the assumptions of Proposition \ref{p6.3}, we have that 
\be
(-1)^{n}\chi(\U)\geq 0.
\ee
\ec
\begin{proof}
Since $\U$ is an affine variety, $\U$ is homotopy equivalent to a finite $n$-dimensional CW complex. Hence $H_{n}(\U^{c})$ is a free $\Gamma$-module. Proposition \ref{p6.3} shows that $H_{i}(\U^{c})$ is a torsion $\Gamma$-module for $i\leq n-1$. So \begin{center}
$\chi(\U)=(-1)^{n} \mbox{ rank} H_{n}(\U^{c}),$
\end{center}
hence $(-1)^{n}\chi(\U)= \mbox{ rank} H_{n}(\U^{c}) \geq 0$.
\end{proof}


\section{Resonance varieties and straightness}\label{S7}

In this section, we introduce the resonance varieties, and describe their relation with characteristic varieties. As before, we assume that $X$ is a finite connected $n$-dimensional CW complex. Consider the cohomology algebra $A^{\bullet}=H^{\bullet}(X,\C)$, with $\dim_{\C} A^{1}=r$. For each $a\in A^{1}$, we have $a^{2}=0$. Then, multiplication by $a$ defines a cochain complex:
\begin{center}
$(A^{\bullet}, \cdot a): A^{0}\overset{a}{\longrightarrow} A^{1} \overset{a}{\longrightarrow} A^{2} \overset{a}{\longrightarrow}\cdots$,
\end{center}
known as the {\it Aomoto complex}.
\bd \label{d7.0} The {\it resonance varieties of $X$} are defined by:\begin{center}
$\re^{i} (X)=\lbrace a\in A^{1} \mid \dim_{\C} H^{i}(A^{\bullet},\cdot a)\neq 0 \rbrace$
\end{center}
\ed
Note that, if $A^{i}=0$, then $\re^{i}(X)=\emptyset$. 
The sets $\re^{i}(X)$ are homogeneous algebraic varieties of the affine space $A^{1}=\C^{r}$ (see \cite[Lemma 3.2]{Su12A} ), and they are homotopy invariants of $X$, see \cite[Lemma 3.3]{Su12A}. Set $$\re(X)= \bigcup_{i=0}^{n}\re^{i}(X).$$


\subsection{Tangent cone inclusion and locally straight space}\label{s7.1}
Assume that $W \subset (\C^{\ast})^{r}$ is a Zariski closed set. Let $J$ be the ideal in the ring of analytic functions $\C \lbrace t_{1},\cdots,t_{r} \rbrace$ defining the germ of $W$ at 1, and let $in(J)$ be the ideal in the polynomial ring $\C[t_{1},\cdots,t_{r}]$ spanned by the initial forms of non-zero elements of $J$.
\bd The {\it tangent cone} of $W$ is defined by \begin{center}
$TC_{1}(W)=V(in(J))$.
\end{center}
\ed

Libgober \cite{L02} established a connection between characteristic varieties and the resonance varieties by the following tangent cone inclusion: 
\bt \label{t7.1} Let $X$ be a finite connected $n$-dimensional CW complex. Then, for any $k\leq n$, \be \label{7.1}
 TC_{1} ( \bigcup _{i=0}^{k} \W_{i}(X)) \subset \bigcup _{i=0}^{k} \re^{i}(X).
 \ee
\et 

\bd  \cite[Definition 6.1]{Su12A} 
Assume that $X$ is a finite connected $n$-dimensional CW complex.
We say $X$ is {\it locally $j$-straight} if the following conditions are satisfied for each $k \leq j$:
\begin{itemize}
\item[(a)] All components of $\bigcup _{i=0}^{k} \W_{i}(X)$ passing through the origin are algebraic subtori.
\item[(b)] $TC_{1} ( \bigcup _{i=0}^{k} \W_{i}(X)) = \bigcup _{i=0}^{k} \re^{i}(X)$.
\end{itemize}
If conditions $(a)$ and $(b)$ hold for all $j\geq 1$, we say $X$ that is a {\it locally straight space}.
\ed

\br\label{BW} The variety $M^{\ast}=\CP \setminus V$ is affine and smooth, so, by \cite[Theorem 1.1]{BW}, the characteristic varieties $\V_{i}(M^{\ast})$ are finite unions of torsion translated subtori. Recall also that by Theorem \ref{t2.1} we have that 
$$\bigcup_{i=0}^{k} \W_{i}(M^{\ast})=( \bigcup_{i=0}^{k} \V_{i}(M^{\ast}) ) \cap \Hom(G^{\ast},\C^{\ast})^{0}.$$
So  condition $(a)$ is always satified for $M^{\ast}$.  Similar considerations apply to the spaces $\U=\CP \setminus (V \cup H)$ and $M=\CN \setminus f^{-1}(0)$, where $f$ is the homogeneous polynomial defining $V$.
\er


\subsection{Divisibility results for resonance varieties}\label{s7.2}
First, we establish a relation between the resonance varieties $\re^{i}(M^{\ast})$ and $\re^{i}(M)$ of $M^*$ and $M$, respectively. Note that $\re^{i}(M) \subseteq H^1(M)\cong \C^r$, while $\re^{i}(M^*) \subseteq H^1(M^*)\cong \C^{r-1}$. However, the Hopf map $p:M \to M^*$ induces a monomorphism $p^{\sharp}:H^1(M^*) \to H^1(M)$, so we can regard $\re^{i}(M^{\ast})$ also as subsets of $H^1(M)$. More precisely, 
$V(\prod_{i=1}^{r} t_{i}^{d_{i}}=1)$  has $\gcd(d_{1},\cdots,d_{r})$ irreducible components in $(\C^*)^r$, with $V(\prod_{i=1}^{r} t_{i}^{d^{\prime}_{i}}=1)$ being the component passing through $(1, \cdots, 1)$. (Recall here that $d^{\prime}_{i}=d_i/\gcd(d_{1},\cdots,d_{r})$.) The linear subspace associated to $V(\prod_{i=1}^{r} t_{i}^{d^{\prime}_{i}}=1)$ in $\C^{r}$ is defined by $\sum_{i=1}^{r} d_{i}^{\prime} z_{i}=0$, and we denote by $L$ this $(r-1)$-dimensional $\C$-vector subspace. Then $L\cong p^{\sharp}(H^1(M^*))$ and, by Propositions \ref{p5.3}, \ref{p2.4}, and the tangent cone inclusion, we have that $$\re(M^{\ast})\subseteq L.$$ If, in addition, $\chi(M^*) \neq 0$, then Proposition \ref{p5.3} and the tangent cone inclusion yield that $\re(M^{\ast})= L.$

\bp \label{p7.1} For any $i\geq 0$, \be  \label{7.2}
\re^{i}(M)=\re^{i}(M^{\ast})\cup \re^{i-1}(M^{\ast}). 
\ee
\ep 
\begin{proof}
We have an isomorphism of graded algebras (see \cite[Proposition 6.4.1]{D2}): \be  \label{7.3}
H^{\bullet}(M) \cong H^{\bullet}(M^{\ast})\otimes H^{\bullet}(\C^{\ast}).
\ee 
Then  
$$\re^{i}(M)=\bigcup_{p+q=i} \re^{p}(M^{\ast})\times \re^{q}(\C^{\ast})= \re^{i}(M^{\ast})\cup \re^{i-1}(M^{\ast}),$$
where the first equality follows from the (proof of the) product formula in \cite[Proposition 13.1]{PS10}, and the last equality follows from the fact that $\re^{q}(\C^{\ast})=\lbrace 0 \rbrace$ for $q=0,1$, and empty, otherwise. 
 \end{proof}
 
 The following corollary is a direct consequence of Proposition \ref{p7.1}.
\bc \label{c7.1}  For any $k\geq 0$, \be \label{7.4}
\bigcup_{i=0}^{k}\re^{i}(M) =\bigcup^{k}_{i=0} \re^{i}(M^{\ast}).
\ee
Moreover, $M^{\ast}$ is a locally $j$-straight space if and only if $M$ is so.
\ec
\begin{proof} The identification in (\ref{7.4}) follows at once from (\ref{7.2}). 
As already pointed out in Remark \ref{BW}, the locally straightness of $M^{\ast}$ and $M$ only depend on condition $(b)$ for each of these spaces. Note that Theorem \ref{t5.1}, combined with Theorem \ref{t2.1}, yields that 
$$\bigcup_{i=0}^{k} \W_{i}(M^{\ast})=( \bigcup_{i=0}^{k} \V_{i}(M^{\ast}) ) \cap \Hom(G^{\ast},\C^{\ast})^{0}= (\bigcup_{i=0}^{k} \W_{i}(M)) \cap \Hom(G^{\ast},\C^{\ast})^{0}.$$
Then the straightness claim follows from the equality (\ref{7.4}).
\end{proof}

In view of the considerations preceding Proposition \ref{p7.1}, Corollary \ref{c7.1} yields that: $$\re(M)=\re(M^{\ast})\subseteq L.$$
 For the local case,  the map $\nu_{x}$ induces an embedding: $\nu_{x}^{\sharp}: \Q^{r_{x}}\hookrightarrow H^{1}(M^{\ast}_{x},\C)$.
Motivated by the equality  $ \bigcup_{i=0}^{k} \W_{i}(M^{\ast}_{x},\nu_{x})   =im(\nu_{x}^{\ast}) \cap (\bigcup_{i=0}^{k} \V_{i}(M^{\ast}_{x}))  $, we now define {\it local resonance varieties} by $$\re(M^{\ast}_{x},\nu_{x})=im(\nu_{x}^{\sharp})\cap \re(M^{\ast}_{x}).$$ The corresponding uniform local resonance variety $\re^{unif}(M^{\ast}_{x},\nu_{x})$ is then defined by
$$\re^{unif}(M^{\ast}_{x},\nu_{x})= \re(M^{\ast}_{x},\nu_{x})\times\C^{r-r_{x}}.$$

\bt  \label{t7.2} Assume that the hypersurface complement $M^{\ast}=\CP \setminus V$ is a locally $j$-straight space. Then, for any $k\leq \min\lbrace j,n-1\rbrace$, we have
\be \label{7.5}
 \bigcup_{i=0}^{k}  \re^{i}(M)= \bigcup_{i=0}^{k}  \re^{i}(M^{\ast}) \subset ( \bigcup_{S\subset V_{1}} \re^{unif}(M^{\ast}_{x_{S}},\nu_{x_{S}}) ) \cap  L,
\ee
where the union is over the strata $S$ of $V_{1}$ with $\dim_{\C} S\geq n-k-1$ and $x_{S}\in S$ is arbitrary.  
\et
\begin{proof} For $k\leq \min\lbrace j,n-1\rbrace$, we have that \[
\begin{aligned}
\bigcup_{i=0}^{k}  \re^{i}(M^{\ast}) 
&  \overset{(1)}{ = } TC_{1} (\bigcup_{i=0}^{k} \W_{i}(M^{\ast})) \\
&  \overset{(2)}{ \subset } TC_{1}\lbrace  (\bigcup_{S\subset V_{1}} \W^{unif}(M^{\ast}_{x_{S}}, \nu_{x_{S}})) \cap V(\prod_{i=1}^{r} t_{i}^{d^{\prime}_{i}}=1) \rbrace  \\
&  \overset{(3)}{ \subset } ( \bigcup_{S\subset V_{1}} \re^{unif}(M^{\ast}_{x_{S}},\nu_{x_{S}}) ) \cap  L 
\end{aligned}
\]
where (1) follows from the locally $j$-straightness of $M^{\ast}$,  (2) follows from (\ref{5.4}) and the fact that, if $W_{1} \subset W_{2}$ then $TC_{1}(W_{1})\subset TC_{1}(W_{2})$, and (3) is a consequence of the tangent cone inclusion.
\end{proof}

 Let $V$ be a hyperplane arrangement in $\CP$, with complement $M^{\ast}$. Then it is known that $M^{\ast}$ is locally straight, see \cite[Proposition 11.1]{Su12A}.  The following result follows from  Theorem \ref{t5.4} by using an argument similar to that of Theorem \ref{t7.2}:
\bt  \label{t7.3} Let $M^{\ast}$ be the complement of an essential hyperplane arrangements. Then, for any $k\leq n-1$,
\be \label{7.9}
 \bigcup_{i=0}^{k}  \re^{i}(M)= \bigcup_{i=0}^{k}  \re^{i}(M^{\ast}) \subset  \bigcup_{S\subset V_{1}} \lbrace \re(M^{\ast}_{x_{S}}) \times V(\sum _{i\notin I_{x_S}} z_{i}=0) \rbrace ,
\ee
where the union is over the strata $S$ of $V_{1}$ with $\dim_{\C} S\geq n-k-1$ and $x_{S}\in S$ is arbitrary.
\et

It is well known that  resonance varieties of hyperplane arrangement  complements are determined  by the intersection lattice. A basic open problem is to find concrete formulas expressing this dependence. Theorem \ref{t7.3} sheds some light into this problem.



\begin{thebibliography}{ADMSP}

\bibitem[A97]{A} D. Arapura, {\it Geometry of cohomology support loci for local systems}, I. J. Algebraic Geom. {\bf 6} (1997), no. 3, 563--597.

\bibitem[B83]{B} A. Borel, {\it Intersection cohomology},
Notes on the seminar held at the University of Bern, Bern, 1983. Reprint of the 1984 edition. Modern Birkh\"{a}user Classics. Birkh\"{a}user Boston, Inc., Boston, MA, 2008.

\bibitem[Bry86]{Br} J.-L. Brylinski, {\it Transformations canoniques, dualit\'e projective, th\'eorie de Lefschetz, transformations de Fourier et sommes trigonom\'etriques},
Ast\'erisque {\bf 140-141} (1986), 3--134, 251.

\bibitem[Bud15A]{Bu} N. Budur, {\it Bernstein-Saito ideals and local systems},
Ann. Inst. Fourier {\bf 65} (2015), no. 2, 549--603. 

\bibitem[Bud15B]{Bu2} N. Budur, {\it Private communication}, 2015.

\bibitem[BW14]{BW} N. Budur, B. Wang,
 {\it Cohomology jump loci of quasi-projective varieties},
 Ann. Sci. \'{E}c Norm. Sup\'{e}r. (4) {\bf 48} (2015), no. 1, 227--236. 
 
 \bibitem[BW15]{BW15} N. Budur, B. Wang,
 {\it Local systems on analytic germ complements}, arXiv:1508.07867.

\bibitem[CDO02]{CDO}  D. Cohen, A. Dimca, P. Orlik, 
{\it Nonresonance conditions for arrangements},
Ann. Inst. Fourier (Grenoble) {\bf 53} (2003), no. 6, 1883--1896.

\bibitem[CS91]{CS} S. E. Cappell and J. L. Shaneson, {\it Singular spaces, characteristic classes, and intersection homology},
Ann. of Math. (2) {\bf 134} (1991), no. 2, 325--374.

\bibitem[DJL07]{DJL}   M. W. Davis, T. Januszkiewicz, I. J. Leary,
{\it The $l^{2}$-cohomology of hyperplane complements}, Groups Geom. Dyn. {\bf 1} (2007), no. 3, 301--309.

\bibitem[Di92]{D1} A. Dimca, {\it Singularities and Topology of Hypersurfaces},
Universitext, Springer-Verlag, New York 1992.

\bibitem[Di04]{D2} A. Dimca, {\it Sheaves in Topology},
Universitext, Springer-Verlag, Berlin, 2004.

\bibitem[Di02]{D3} A. Dimca, 
{\it Hyperplane arrangements, M-tame polynomials and twisted cohomology}, in: 
Commutative algebra, singularities and computer algebra (Sinaia, 2002), 113--126, NATO Sci. Ser. II Math. Phys. Chem. {\bf 115}, Kluwer Acad. Publ., Dordrecht, 2003.

\bibitem[DL07]{DL} A. Dimca, A. Libgober,
{\it Local topology of reducible divisors},
in: Real and complex singularities, 99--111, Trends Math., Birkh\"{a}user, Basel, 2007.

\bibitem[DM07]{DM} A. Dimca,  L. Maxim,
{\it Multivariable Alexander invariants of hypersurface complements},
Trans. Amer. Math. Soc. {\bf 359} (2007), no. 7, 3505--3528.

\bibitem[DN04]{DN} A. Dimca, A. N\'{e}methi, 
{\it Hypersurface complements, Alexander modules and monodromy}, in: Real and complex singularities, 19--43, Contemp. Math. {\bf 354}, Amer. Math. Soc., Providence, RI, 2004.

\bibitem[DP03]{DP} A. Dimca, S. Papadima,
{\it Hypersurface complements, Milnor fibers and higher homotopy groups of arrangments},
Ann. of Math. (2) {\bf 158} (2003), no. 2, 473--507.

\bibitem[DF87]{DF} W. G. Dwyer, D. Fried,
{\it Homology of free abelian covers. I},
Bull. London Math. Soc. {\bf 19} (1987), no. 4, 350--352.

\bibitem[Huh12]{Huh} J. Huh,
{\it Milnor numbers of projective hypersurfaces with isolated singularities}, Duke Math. J. {\bf 163} (2014), 1525--1548.

\bibitem[Lib82]{L82} A. Libgober, 
{\it Alexander polynomial of plane algebraic curves and cyclic multiple planes},
Duke Math. J. {\bf 49} (1982), no. 4, 833--851.

\bibitem[Lib83]{L83} A. Libgober, 
{\it Alexander invariants of plane algebraic curves},
in: Singularities, Part 2 (Arcata, Calif., 1981), 135--143, Proc. Sympos. Pure Math., {\bf 40}, Amer. Math. Soc., Providence, RI, 1983.

\bibitem[Lib92]{L92} A. Libgober, 
{\it On the homology of finite abelian coverings}, Topology Appl. {\bf 43} (1992), no. 2, 157--166.

\bibitem[Lib94]{L94} A. Libgober,
 {\it Homotopy groups of the complements to singular hypersurfaces, $\uppercase\expandafter{\romannumeral 2}$},
Ann. of Math. (2) {\bf 139} (1994), no. 1, 117--144.

\bibitem[Lib01]{L01} A. Libgober,
 {\it Characteristic varieties of algebraic curves},
in: Applications of algebraic geometry to coding theory, physics and computation (Eilat, 2001), 215--254, NATO Sci. Ser. II Math. Phys. Chem. {\bf  36}, Kluwer Acad. Publ., Dordrecht, 2001. 

\bibitem[Lib02]{L02} A. Libgober,
 {\it First order deformations for rank one local systems with a non-vanishing cohomology},
in: Arrangements in Boston: a Conference on Hyperplane Arrangements (1999). Topology Appl. {\bf 118} (2002), no. 1-2, 159--168. 

\bibitem[Lib09]{L09} A. Libgober,
 {\it Non vanishing loci of Hodge numbers of local systems},
 Manuscripta Math. {\bf 128} (2009), no. 1, 1--31.

\bibitem[Liu16]{Liu} Y. Liu, {\it Nearby cycles and Alexander modules of hypersurface complements}, Adv. Math. {\bf 291} (2016), 330--361. 

\bibitem[LiM15]{LiM} Y. Liu, L. Maxim, {\it Reidemeister torsion, peripheral complex, and Alexander polynomials of hypersurface complements}, Algebr. Geom. Topol. {\bf 15} (2015), no. 5, 2757--2787.  

\bibitem[Max06]{Max} L. Maxim, 
{\it Intersection homology and Alexander modules of hypersurface complements},
Comment. Math. Helv. {\bf 81} (2006), no. 1, 123--155.

\bibitem[PS10]{PS10} S. Papadima, A. Suciu,
{\it Bieri-Neumann-Strebel-Renz invariants and homology jumping loci}, Proc. Lond. Math. Soc. (3) {\bf 100} (2010), no. 3, 795--834.

\bibitem[Sab90]{Sab} C. Sabbah, {\it Modules d'Alexander et $\mathcal{D}$-modules}, Duke Math. J. {\bf 60} (1990), no. 3, 729--814.

\bibitem[Sch03]{Sc} J. Sch\"{u}rmann, 
{\it Topology of Singular Spaces and Constructible Sheaves},
Monografie Matematyczne {\bf 63}, Birkh\"auser Verlag, Basel 2003.

\bibitem[Ser65]{Ser} J. P. Serre, 
{\it Alg\`ebre locale. Multiplicit\'es},
Lecture Notes in Mathematics {\bf 11}, Springer-Verlag, Berlin-New York 1965.

\bibitem[Su00]{Su00} A. Suciu,
{\it Fundamental groups of line arrangements: enumerative aspects}, in: 
Advances in algebraic geometry motivated by physics (Lowell, MA, 2000), 43--79, Contemp. Math. {\bf 276}, Amer. Math. Soc., Providence, RI, 2001.

\bibitem[Su12A]{Su12A} A. Suciu,
{\it Resonance varieties and Dwyer-Fried invariants}, in: 
Arrangements of hyperplanes-Sapporo 2009, 359--398, Adv. Stud. Pure Math. {\bf 62}, Math. Soc. Japan, Tokyo, 2012. 

\bibitem[Su12B]{Su12B} A. Suciu,
{\it Geometric and homological finiteness in free abelian covers}, in: 
Configuration Spaces: Geometry, Combinatorics and Topology (Centro De Giorgi, 2010), 461--501, Publications of the Scuola Normale Superiore, vol. {\bf 14}, Edizioni della Normale, Pisa, 2012.   

\bibitem[Su14]{Su14} A. Suciu,
{\it Characteristic varieties and Betti numbers of free abelian covers}, 
Int. Math. Res. Not. IMRN {\bf 2014}, no. 4, 1063--1124.

\bibitem[Wei94]{Wei} C. A. Weibel, 
{\it An introduction to homological algebra},
Cambridge Studies in Advanced Mathematics {\bf 38}. Cambridge University Press, Cambridge, 1994. 

\end{thebibliography}
\end{document}